\documentclass[12pt]{article}
\usepackage{amsmath,amsthm,amsfonts,amssymb,amscd}
\def\qed{{\unskip\nobreak\hfil\penalty50
\hskip2em\hbox{}\nobreak\hfil$\square$
\parfillskip=0pt \finalhyphendemerits=0\par}\medskip}
\def\proof{\trivlist \item[\hskip \labelsep{\bf Proof\ }]}
\def\endproof{\null\hfill\qed\endtrivlist}

\def\Ad{{\mathrm {Ad}}}

\def\End{{\mathrm {End}}}
\def\Hom{{\mathrm {Hom}}}

\def\o{{\mathrm {opp}}}

\def\Vir{{\mathrm {Vir}}}

\def\dim{{\mathrm {dim}}}
\def\a{\alpha}
\def\b{\beta}

\def\e{\varepsilon}
\def\f{\varphi}
\def\g{\gamma}

\def\l{\lambda}

\def\phi{\varphi}

\def\th{\theta}

\def\o{\omega}
\def\Om{\Omega}

\def\s{{\sigma}}

\def\r{{\rho}}
\def\t{{\tau}}
\def\z{{\zeta}}

\newtheorem{theorem}{Theorem}[section]
\newtheorem{lemma}[theorem]{Lemma}

\newtheorem{corollary}[theorem]{Corollary}

\newtheorem{proposition}[theorem]{Proposition}

\def\emptyset{\varnothing}
\def\setminus{\smallsetminus}
\def\Vect{{\mathrm {Vect}}}
\def\Diff{{\mathrm {Diff}}}
\def\Mob{{\rm\textsf{M\"ob}}}
\def\PSL{PSU(1,1)}

\def\res{\!\restriction\!}
\def\A{{\cal A}}
\def\B{{\cal B}}
\def\C{{\cal C}}
\def\D{{\cal D}}
\def\O{{\cal O}}

\def\I{{\cal I}}

\def\N{{\cal N}}
\def\M{{\cal M}}

\def\H{{\cal H}}

\def\U{{\cal U}}
\def\V{{\cal V}}

\title{\huge Topological Sectors and a Dichotomy\\
in Conformal Field Theory}
\author{
{\sc Roberto Longo}\footnote{Supported in part by GNAMPA-INDAM and 
MIUR.}\\
Dipartimento di Matematica\\
Universit\`a di Roma ``Tor Vergata''\\
Via della Ricerca Scientifica, 1, I-00133 Roma, Italy\\
E-mail: {\tt longo@mat.uniroma2.it}\\
{}\\
{\sc Feng Xu}\footnote{Supported in part by NSF.}\\
Department of Mathematics\\
University of California at Riverside\\
Riverside, CA 92521\\
E-mail: {\tt xufeng@math.ucr.edu}}
\begin{document}
\date{}
\maketitle

\begin{abstract}
Let $\A$ be a local conformal net of factors on $S^1$ with the split
property.  We provide a topological construction of soliton
representations of the $n$-fold tensor product
$\A\otimes\cdots\otimes\A$, that restrict to true representations of
the cyclic orbifold $(\A\otimes\cdots\otimes\A)^{\mathbb Z_n}$.  We
prove a quantum index theorem for our sectors relating the Jones index
to a topological degree.  Then $\A$ is not completely rational iff the
symmetrized tensor product $(\A\otimes\A)^{\rm flip}$ has an
irreducible representation with infinite index.  This implies the
following dichotomy: if all irreducible sectors of $\A$ have a conjugate 
sector then either $\A$ is completely rational or $\A$ has uncountably many
different irreducible sectors.  Thus $\A$ is rational iff $\A$ is
completely rational.  In particular, if the $\mu$-index of $\A$ is
finite then $\A$ turns out to be strongly additive.  By \cite{KLM},
if $\A$ is rational then the tensor category of representations of
$\A$ is automatically modular, namely the braiding symmetry is
non-degenerate.  In interesting cases, we compute the fusion rules of
the topological solitons and show that they determine all twisted
sectors of the cyclic orbifold.
\end{abstract}

\newpage

\section{Introduction}
The main theme of this paper, topological sectors in Conformal Quantum 
Field Theory, has been the subject of interest by the authors for different 
reasons. 

One motivation came from the study of the sector structure in 
the cyclic orbifold associated with rational models, where the 
operator algebraic methods go beyond the analysis by 
infinite Lie algebras, in particular by using the structure results 
in \cite{KLM}. As we shall see, a quantum index theorem by the Jones index 
captures an essential part of information here.

Another motivation came in relation to irrational Conformal Field 
Theory, where most of the underlying structure is still to be 
uncovered.  Also in this case, the algebraic approach is essential and 
leads to a surprising finite/uncountable dichotomy concerning the 
set of all irreducible sectors in the rational/irrational case.

Before stating further results and consequences of our work, and 
explaining in more detail the above mentioned issues, we recall the 
notion of complete rationality \cite{KLM} which is at the basis of our 
analysis.

In all the present paper we shall deal with \emph{diffeomorphism 
covariant} (irreducible) local nets of von Neumann algebras on $S^1$, 
called conformal nets, and we 
explain our results in this framework, although weaker assumptions 
would be sufficient.

{\it Complete rationality.} 
Let then $\A$ be a local conformal net on $S^1$. 
$\A$ is called completely rational if 
\begin{itemize}
\item $\A$ is split,
\item $\A$ is strongly additive,
\item The $\mu$-index $\mu_{\A}$ is finite.
\end{itemize}
The first two conditions are, in a certain sense, one another dual. If $I_1,I_2$ 
are intervals of $S^1$, the split property states that the local von 
Neumann algebras $\A(I_1)$ and $\A(I_2)$ ``maximally decouple'' if $I_1$ and $I_2$ 
have disjoint closures, namely  $\A(I_1)\vee\A(I_2)$ is naturally 
isomorphic to $\A(I_1)\otimes\A(I_2)$;
while strong additivity requires that $\A(I_1)$ and 
$\A(I_2)$ ``maximally interact'' if $I_1$ and $I_2$ have a common 
boundary point, namely $\A(I_1)\vee\A(I_2)=\A(I)$ where $I$ is the 
union of $I_1$, $I_2$ and the boundary point, see e.g. \cite{L4} and 
refs. therein \footnote{We shall later see that, 
in the diffeomorphism covariant case, strong additivity follows from 
the other two conditions.\par 
The symbol ``$\vee$'' denotes the von Neumann algebra generated.}.

In the last condition, $\mu_{\A}$ is the 
Jones index \cite{J} of the inclusion of factors $\A(E)\subset \A(E')'$ were 
$E\subset S^1$ and its complement $E'$ are union of two proper 
disjoint intervals.

In \cite{Xu6} it was shown that $\mu_{\A}<\infty$ when $\A$ is 
associated with $SU(N)$ loop group models. The general theory of 
complete rationality was developed in \cite{KLM}. To check complete 
rationality one may use the fact that this property equivalently 
holds for finite-index subnet \cite{L4}.

One of the main points is that, if $\A$ is completely rational, then
\begin{equation}\label{sum}
\mu_{\A}=\sum_i d(\rho_i)^2\ ,
\end{equation}
the $\mu$-index equals the global index i.e.
the sum of the indeces (= squares of dimensions)
of all irreducible sectors; thus $\A$ is 
rational and indeed the representation tensor category is even modular.

One issue in this paper is to extend the above equality to non rational 
nets. This will lead in particular to a general characterization of 
rational nets.

A look at the basic models constructed by positive-energy representations of the 
diffeomorphism group, the Virasoro nets $\Vir_c$, gives insight here. If the 
central charge $c$ is less then one then $\Vir_c$ is completely 
rational, as is indeed the case of all conformal nets with $c<1$ \cite{KL}. By 
contrast, if $c>1$ then $\Vir_c$ is not even strongly additive \cite{BS} and 
has uncountably many sectors as is known, see e.g. \cite{C}. 
The boundary case, $\Vir_c$ at $c=1$, has uncountably many 
sectors and has recently 
been shown in  \cite{Xu4} to be strongly additive. 
Moreover in the case $c\geq 1$ 
there are plenty of infinite index sectors \cite{C}. We shall see that 
the structure manifested by Virasoro nets undergoes a general phenomenon.
But, before this, we give a general picture of our mentioned dichotomy.

{\it The dichotomy.} Dichotomies concerning the cardinality of various 
structures appear in Mathematics.  One simple example concerns a
$\sigma$-algebra: it is either finite or uncountable.  This is an 
immediate consequence of the basic Cantor-Bernstein theorem to the 
effect that $2^{\mathbb N}$ is uncountable.

Also elementary is the statement that a Hamel basis for a 
Banach space is either finite or uncountable. This is due to 
Baire category theorem. The dichotomy holds because limit points 
are included in the structure.

One further example is provided by a compact group. Again it is 
either finite or uncountable. Here the statement follows at once by the 
existence of a finite Haar measure, a structure property of global nature.

As a final example consider the case of a separable, simple, unital
$C^*$-algebra $\mathfrak A$ and denote by Irr$\mathfrak A$ the set of
equivalence classes of irreducible representations of $\mathfrak A$. 
Then either Irr$\mathfrak A$ consists of a single element ($\mathfrak A$ is a
matrix algebra) or Irr$\mathfrak A$ is uncountable.  This fact is a
consequence of the deep theorem of Glimm on the classification of type
$I$ $C^*$-algebras \cite{Dix} (if $\mathfrak A$ has a representation
not of type $I$ then uncountably many irreducibles have to appear in
its disintegration).

The dichotomy in this paper is more similar in the spirit to this last 
example: a high degree of understanding of the structure is necessary to get it.

The statement is the following.
Let $\A$ be a local conformal net with the split property. Assume 
that every irreducible sector of $\A$ has a conjugate sector. Then either 
$\A$ is completely rational or $\A$ admits uncountably many different 
irreducible sectors. We shall later return on the consequences of 
this fact.

Now, to exhibit uncountably many sectors in the irrational 
case, some new construction of representations has to appear at some 
stage. This is indeed one of the most interesting points. 
These representations are constructed topologically, as we now explain.
(Note that, in higher dimension spacetimes, charges of topological 
nature and wide localization have long been known and are natural 
in particular in quantum electrodynamics, see \cite{BF}).

{\it Topological sectors.} Let's start with a simple observation. 
Let $\A$ be a conformal net and $\A_0$ its 
restriction to the real line $\mathbb R = S^1\setminus\{\z\}$ 
obtained by removing a point 
$\z$ from the circle. If $h:\mathbb R\to S^1$ is a smooth, injective, 
positively oriented map, we get a representation $\Phi_h$ of the $C^*$-algebra 
$\overline{\cup_I \A_0(I)}$ (union over all bounded intervals of 
$\mathbb R$) by setting
\[
\Phi_h(x)\equiv U(k_I)xU(k_I)^*,\quad x\in\A_0(I),
\]
where $k_I:S^1\to S^1$ is any diffeomorphism of $S^1$ that 
coincides with $h$ on the interval $I$, 
and $U$ is the covariance projective unitary 
representation of $\Diff(S^1)$. Assuming $h$ to be smooth also at $\pm 
\infty$, then $\Phi_h$ is a \emph{soliton}, namely it is normal on the 
algebras associated with half-lines. 

Incidentally, this gives an elementary and model 
independent construction of type III representations, see \cite{DS,CR} 
for constructions of type III representations in models.

Let now $f:S^1\to S^1$ be a smooth, locally injective map of degree 
$n = \text{deg}f \geq 1$. Then $f$ has exactly $n$ right inverses 
$h_i$, $i=0,1,\dots n-1$, namely there are $n$ injective 
smooth maps $h_i:S^1\setminus\{\z\}\to S^1$ such that $f(h_i(z))=z$, 
$z\in S^1\setminus\{\z\}$. The $h_i$'s are smooth also at $\pm\infty$.
For the moment we make an arbitrary choice of order $h_0, h_1, h_2,\dots $. 

As just explained, we have $n$ solitonic 
representations $\Phi_{h_i}$ of $\A$, hence one (reducible) soliton
$\Phi_f\equiv\Phi_{h_0}\otimes\cdots\otimes\Phi_{h_{n-1}}$ of $\A\otimes\cdots\otimes\A$.

Now, if $I\subset\mathbb R$ is an interval, the intervals $I_i\equiv 
h_i(I)\subset S^1$ have pairwise disjoint closures hence, by the split 
property that we now assume, there is a natural identification 
\[
\chi_I: \A(I_0)\otimes\cdots\otimes\A(I_{n-1})\to 
\A(I_0)\vee\cdots\vee\A(I_{n-1})
\]
therefore we get an irreducible solitonic representation 
$\pi_f$ of $\A\otimes\cdots\otimes\A$ by gluing together the range of $\Phi_{f,I}$ 
by $\chi_I$, namely 
\[
\pi_{f,I}\equiv\chi_I\cdot\Phi_{f,I}\ .
\]
Let's say now that we choose the $h_i$'s so that 
the sequence of intervals $I_0, I_1,\dots I_{n-1}$ is 
counter-clockwise increasing (this requirement does not depend on $I$).
This fixes the order of the $h_i$'s up to a cyclic permutation.

If we go to the cyclic orbifold $(\A\otimes\cdots\otimes\A)^{\mathbb Z_n}$
the dependence on the cyclic permutations disappears and we can 
easily verify that
\[
\t_f \equiv \pi_f\restriction (\A\otimes\cdots\otimes\A)^{\mathbb Z_n}
\]
is indeed a well-defined \emph{DHR representation} with exactly $n$ 
irreducible components. We have thus generated a 
family of twisted sectors for the cyclic orbifold.

We shall further see that $\pi_f$ itself does not depend, up to 
unitary equivalence, on the ordering of the 
$h_i$'s, by choosing the $I_i$'s increasing as 
above, namely the soliton sector $[\pi_f]$ is an intrinsic object.
In other words, if we denote by
$\pi_{f,p}$ the sector corresponding to another other ordering, 
where $p\in\mathbb P_n$ is the permutation rearranging the $h_i$'s, 
then $[\pi_{f,p}]$ depends only on the cosets $\mathbb 
P_n/\mathbb Z_n$. The conjugate sector of $\pi_f$ corresponds to the 
clockwise ordering of the $h_i$'s.

{\it A quantum index theorem.}  The soliton representation $\pi_f$, 
and its DHR restriction $\t_f$,
depend on $f$ only up to unitary equivalence. In a sense 
these topological sectors play a role similar to the Toeplitz 
operators (see e.g. \cite{A}) in the framework of Fredholm linear
operators, where the analytical index coincides with the degree deg$f$. 
\footnote{Here, however, the underlying space (net) depends on deg$f$.}
As explained in \cite{L6}, Doplicher-Haag-Roberts 
localized endomorphisms \cite{DHR} may be viewed as a Quantum Field Theory 
analog of elliptic operators, in the context of a quantum index theorem. 
The topological sectors provide a good illustration of this point. 
Denoting by $\t^{(i)}_f$, $i\in\mathbb Z_n$, the $n$ direct summands 
of $\t_f$, we have
\begin{equation*}
\text{Index}(\t^{(i)}_f) = \text{Index}(\pi_f) = \mu_\A^{n-1}\ .
\end{equation*}
Here the index is the Jones index \cite{J}, the analog of the Fredholm 
index \cite{L1,L1'}, $\mu_\A$ is the above structure constant for $\A$, 
and $n = \text{deg}f$ is the topological index, which is manifestly deformation 
invariant. 

As we shall see, more general topological sectors arise from non-vacuum 
representations. The index and further structure of these sectors will 
be determined.

Most of the results in this paper depend, maybe implicitly, on the 
above index formulas.

{\it The structure of the sectors.} At the infinitesimal level, the 
twisted sectors of the cyclic orbifold have already been considered in the 
papers \cite{BDM, BHS} in the framework of Kac-Moody Lie algebras. To 
study the structure of the tensor category of topological sectors, 
it is however necessary to have the sectors in global exponential form 
and a general theory at one's disposal, as provided by our approach. 

In Section \ref{fus} we shall determine all the twisted 
irreducible sectors of the $n$-cyclic orbifold and 
give a detailed account of the 
fusion rules of the topological solitons, in the cases $n=2,3,4$, for a 
general completely rational net. The method of $\a$-induction 
\cite{LR,Xu1,BE,BEK} is here essential.

What undergoes the structure of sectors is the covariance symmetry. The 
Lie algebra of $\Diff(S^1)$, namely the Virasoro algebra at $c=0$ with 
generators $L_n$ and relations 
\[
[L_m,L_n]=(m-n) L_{m+n}\ ,
\]
has an endomorphism, for each positive integer $k$, given by 
\[
L_n \mapsto\frac{1}{k}L_{kn}\ ,
\]
see \cite{BHS,SW}. As we shall see, this 
corresponds to an embedding of the $k$-cover of $\Diff(S^1)$ 
into $\Diff(S^1)$. The covariance projective unitary representations 
are obtained by composing with this embedding the original 
representation (of the vacuum or of a non-trivial sector).

{\it Rationality, modularity, strong additivity, sectors with infinite index.}
The above described dichotomy has the following corollary.  A 
conformal net $\A$ with the split property is rational, in the sense 
that the representation tensor category has only finitely many inequivalent 
irreducible objects and all have a conjugate, if and only if $\A$ is 
completely rational.

All results obtained for completely rational nets \cite{KLM} then immediately 
apply to rational conformal nets. Among them we only mention that a 
rational conformal net has a modular representation tensor category, 
namely the braiding symmetry is automatically non-degenerate. The 
modularity property is at the basis of most of the analysis in 
Conformal QFT and is often taken for granted or implicitly 
conjectured to hold, see e.g. \cite{GW}.

Note that our work shows in particular that a conformal net with the
split property and finite $\mu$-index is automatically strongly
additive or, equivalently, Haag dual on the real line.  To understand
the interest of this result, note that the strong additivity property
is crucial for many results and often one of the hardest point to
prove, see \cite{Xu3,L4,Xu4,C}.  As suggested by Y. Kawahigashi, 
strong additivity can be thought as an amenability 
property; our result supports this view. Our proof makes use of basic
properties of simple subfactors \cite{L0}.

A further consequence is that the equality (\ref{sum}) 
between the $\mu$-index and the global index holds true
for any diffeomorphism covariant local net 
with the split property (regardless $\mu_\A$ is finite or infinite), a 
non-trivial useful result at the basis of our work. 

Last, we state the following characterization of being not
complete rationality: a conformal net $\A$ with the split property is 
not completely rational if and only if the $2$-orbifold net
$(\A\otimes\A)^{\text{flip}}$ admits an irreducible sector with 
infinite index. 

General properties of sectors with infinite index were studied in 
\cite{GL1}, but first examples were constructed by Fredenhagen in 
\cite{F}. Indeed, as mentioned, Carpi \cite{C1,C} has recently shown that 
irreducible sectors with infinite index appear in the Virasoro nets 
$\Vir_c$ if $c\geq 1$, as suggested by Rehren in \cite{R1}. (By 
contrast notice that, in QFT on Minkowski spacetime, all irreducible DHR 
sectors  with an isolated mass shell have finite dimension  \cite{BF}).

Our general construction of infinite index irreducible sectors is 
natural and surprising.  Consider indeed the case of the net $\A$ 
associated with the $U(1)$-current algebra.  All sectors are known in 
this case \cite{BMT}, the irreducible ones all have index 1 and form a 
one-parameter family.  Thus, by \cite{KLM}, $\A\otimes\A$ has only a 
two-parameter family of irreducible sectors, all with index 1.  Yet, 
the ``trivial'' passage to the index 2 subnet 
$(\A\otimes\A)^{\text{flip}}$ makes infinite index sectors to appear.  
Note that one of Fredenhagen's examples is similar in the spirit, but concerns 
the infinite index subnet $(\A\otimes\A)^{SO(2)}$.

At this point we close our expository part and refer to the 
rest of the paper for a detailed account and further results.
See  \cite{H,T} as reference books.
\section{On the symmetry groups}
\label{first}
We shall denote by $\Diff(S^1)$ the group of orientation preserving 
smooth diffeomorphisms of $S^1\equiv\{z\in\mathbb C: |z|=1\}$.  
$\Diff(S^1)$ is an infinite dimensional Lie
group whose Lie algebra is $\Vect(S^1)$, the Lie algebra of smooth 
vector fields on the circle. The complexification $\Vect_{\mathbb 
C}(S^1)$ of $\Vect(S^1)$ has a basis with elements 
$L_n\equiv -z^{n+1}\frac{\rm d}{{\rm d}z}$, $n\in\mathbb Z$, 
satisfying the relations
\begin{equation}\label{vectrel}
[L_m,L_n]=(m-n) L_{m+n} . 
\end{equation}
We shall consider $\Vect_{\mathbb C}(S^1)$ as a Lie algebra with 
involution
\[
L_n^* = L_{-n} .
\]
$SU(1,1)$ is the group of $2\times 2$ matrices defined 
by:
\begin{equation}
 SU(1,1)\equiv\Big \{ \begin{pmatrix}\alpha & \beta \\
 \bar{\beta} & \bar{\alpha}
 \end{pmatrix} \Big|
 \:\alpha,\beta \in \mathbb C,\; 
|\alpha|^{2}-|\beta|^{2}=1\Big \}.
\end{equation}
$SU(1,1)$ acts on $S^1$ by linear 
fractional transformations:
\begin{equation}
g(z)\equiv \frac{\a z + \b}{\bar\b z +\bar\a}\ ,
\end{equation}
where $g=\begin{pmatrix}\alpha & \beta \\
 \bar{\beta} & \bar{\alpha}\end{pmatrix}$.
This action factors through a faithful action of $PSU(1,1)\equiv 
SU(1,1)/\{ \pm1 \}$ on $S^1$.
 
The corresponding diffeomorphisms $z\mapsto g(z)$ constitute a 
subgroup of 
$\Diff(S^1)$, the \emph{M\"{o}bius group} 
$\Mob$, isomorphic to $PSU(1,1)$.
$SU(1,1)$ is a double cover of $PSU(1,1)$, thus
$PSU(1,1)$ and $SU(1,1)$ are locally isomorphic, they have the same 
Lie algebra 
$s\ell (2,\mathbb R)$. 
The complexified Lie algebra 
$s\ell(2,\mathbb C)$ of 
$s\ell(2,\mathbb R)$ has generators $L_{-1},L_{0},L_{1}$ satisfying 
the relations
\[
[L_{1},L_{-1}]=2L_{0}, \qquad [L_{\pm 1},L_{0}]=\pm L_{\pm 1}.
\]
Therefore the elements $L_{-1},L_{0},L_{1}$ of $\Vect_{\mathbb 
C}(S^1)$ 
exponentiate to a subgroup of $\Diff(S^1)$ locally isomorphic to 
$\Mob$. As exp$(2\pi L_0)$ is the identity of $\Diff(S^1)$, this 
group is 
indeed isomorphic to $\Mob$.
$\Vect_{\mathbb C}(S^1)$  contains infinitely many
further copies of $s\ell(2,\mathbb C)$; for a fixed integer $n>0$ we 
get a copy generated by the elements
$L_{-n},L_{0},L_{n}$. Setting 
\begin{equation}
  L'_{m} \equiv \frac{1}{|n|}L_{m},\quad m = n, -n, 0, 
\end{equation}         
we have indeed the relations
\begin{equation}
  [L'_{n},L'_{-n}]=2L'_{0},\qquad
         [L'_{\pm n},L'_{0}]= \pm L'_{\pm n} .
\end{equation}
The Lie subgroup of $\Diff(S^{1})$ corresponding to 
$L_{-n},L_{0},L_{n}$ is thus a cover of $\Mob$. As $\text{exp}(2\pi 
nL'_0)=\text{exp}(2\pi L_0)$ 
is the identity, this group is then isomorphic to $\Mob^{(n)}$,  
the $n$-cover group of $\Mob$. Thus there is a natural embedding
\[
M^{(n)}: \Mob^{(n)} \hookrightarrow \Diff(S^1)\ .
\]
With $g=\left(\begin{smallmatrix}\alpha & 
\beta\\
 \bar{\beta} & \bar{\alpha}
 \end{smallmatrix}\right) \in SU(1,1)$ we shall see that 
 $M^{(n)}_g\in \Diff(S^1)$ is formally given by
\begin{equation}\label{nsqrt}
  M^{(n)}_{g}(z)\equiv \sqrt[n]{\frac{\alpha z^{n} + \beta}
             {\bar{\beta}z^{n}+\bar{\alpha}}}\ ,
\end{equation}
if we locally identify $SU(1,1)$ and $\Mob^{(n)}$.

Denote by $g\mapsto\underline  g$ the quotient map 
$\Mob^{(n)}\to\Mob$.
\begin{proposition}\label{cover}
There is a unique continuous isomorphism $M^{(n)}$ of \Mob$^{(n)}$ 
into 
$\Diff(S^1)$ such that the following diagram commutes for every 
$g\in\text{\Mob}^{(n)}$
\begin{equation}\label{cdn}
\CD
S^1 @>M^{(n)}_g>> S^1   \\
@V z^n V V @V V z^n V\\
S^1 @>M_{\underline g}>> S^1
\endCD
\end{equation}
i.e. $M_g^{(n)}(z)^n=M_{\underline g}(z^n)$ for all $z\in S^1$.
\end{proposition}
Denote by $\sqrt[n]{z}$ the $n^\text{th}$-root function 
on the cut plane $\mathbb C\setminus (-\infty,0]$.
For a fixed $g\in PSU(1,1)$, the map 
$f_g: z\in S^1\mapsto\frac{\alpha z^{n} + \beta}
{\bar{\beta}z^{n}+\bar{\alpha}}\in S^1$ has winding number $n$. 
The Riemann surface $\Sigma_n$ associated with the function 
$\sqrt[n]{z}$ is a 
$n$-cover of $\mathbb C\setminus\{0\}$, we may thus lift $f_g$ to a 
one-to-one  map 
$\tilde{f}_g$  from $S^1$ to the elements of $\Sigma_n$ projecting 
onto $S^1$ on $\mathbb C\setminus\{0\}$. The lift is uniquely 
determined 
as soon as we specify the value $\tilde f_g(1)$ among the $n$ 
elements of $\Sigma_n$ projecting onto $f_g(1)$. 

Let $\V$ be a connected neighborhood of the identity in $PSU(1,1)$ 
such that 
$f_g(1)\in S^{1}\setminus\{-1\}$ for  all $g\in\V$. Then we define 
$\tilde f_g$ for $g\in\V$ by requiring that 
$\tilde f_g(1)=f_{g}(1)\in
\mathbb C\setminus (-\infty,0]\subset\Sigma_n$.

We then set
\begin{equation}\label{Mn}
  M^{(n)}_{g}(z)\equiv \sqrt[n]{\tilde{f}_g(z)},\quad g\in\V.
\end{equation}
Choosing a neighborhood $\V_0$ of the identity in $PSU(1,1)$ such 
that
$\V_0\cdot\V_0\subset\V$ we then have
\[
M^{(n)}_{gh}=M^{(n)}_{g}M^{(n)}_{h}, \quad g,h\in\V_0\ ,
\]
namely we have a local isomorphism of $\V_0\subset PSU(1,1)$ into 
$\Diff(S^1)$, and this extends to a global isomorphism of 
$\Mob^{(n)}$ into 
$\Diff(S^1)$, still denoted by $M^{(n)}$. 

Clearly $M_g^{(n)}(z)^n=M_{\underline g}(z^n)$ for all $z\in 
S^1$ if $g\in\V$.

Note that for any ${\underline g}\in PSU(1,1)$ and $g\in \Mob^{(n)}$ 
projecting onto 
${\underline g}$, we have 
$n$ diffeomorphisms 
\begin{equation}\label{multivocal}
R(\tfrac{2k\pi}{n})M^{(n)}_{g},\quad k=0,1,\cdots n-1\ .
\end{equation}
corresponding to the other possible choices of $g$. Here $R$ is the 
rotation one-parameter subgroup of $\Mob$. Thus $M_g^{(n)}(z)^n$ is 
independent of the pre-image $g$ of $\underline g$. It follows by the 
multiplicative property that $M_g^{(n)}(z)^n=M_{\underline g}(z^n)$ 
for all $z\in 
S^1$ and all $g\in \Mob^{(n)}$.

Concerning the uniqueness of $M^{(n)}$ note that $M_I^{(n)}=I$ 
because 
$M^{(n)}$ is an isomorphism. By continuity $M_g^{(n)}(z)\in 
S^1\setminus \{-1\}$ for $g$ in a neighborhood $\U$ of $I$ and this 
determines  $M^{(n)}$ on $\U$, hence on all $\Mob^{(n)}$.
\endproof

Of particular interest is the case $n=2$. $\Mob^{(2)}$ is isomorphic 
to $SU(1,1)$ and we thus have an isomorphism
\begin{equation}\label{M2}
M^{(2)}: SU(1,1)\hookrightarrow\Diff(S^1)\ .
\end{equation}
We shall often identify $PSU(1,1)$ with $\Mob$ and $SU(1,1)$ with 
$\Mob^{(2)}$. 

We now extend the above proposition to general diffeomorphisms.

Denote by $\Diff^{(n)}(S^1)$ the $n$-central 
cover group of $\Diff(S^1)$. The group $\Diff^{(n)}(S^1)$ is obtained 
from  $\Diff^{(n)}(S^1)$ similarly as $\Mob^{(n)}$ is obtained from 
$\Mob$ 
(the $1$-torus rotation subgroup lifts to its $n$-cover), but we 
shall 
soon give an explicit realization of $\Diff^{(n)}(S^1)$.

The universal cover
group  $\Diff^{(\infty)}(S^1)$ of $\Diff(S^1)$ is the projective limit
\[
\Diff^{(\infty)}(S^1)\equiv \varprojlim_{n\in\mathbb 
N}\Diff^{(n)}(S^1).
\]

If $n\in\mathbb N$, the map
\[
\Vect_{\mathbb C}(S^1) \to \Vect_{\mathbb C}(S^1) ,\quad L_m\mapsto 
\frac{1}{n}L_{nm}
\]
defines an injective endomorphism of $\Vect_{\mathbb C}(S^1)$. Its 
inverse corresponds to a an embedding 
\[
M^{(n)}:\Diff^{(n)}(S^1)\hookrightarrow\Diff(S^1)
\]
that extends the one in (\ref{cdn}) (still denoted by the same 
symbol).

Denote by $g\mapsto\underline  g$ also the quotient map 
$\Diff^{(n)}(S^1)\to\Diff(S^1)$. We then have:
\begin{proposition}\label{cover2}
There is a unique continuous isomorphism $M^{(n)}$ of 
$\Diff^{(n)}(S^1)$ into 
$\Diff(S^1)$ such that the diagram (\ref{cdn}) commutes for every 
$g\in\Diff^{(n)}(S^1)$, namely $M_g^{(n)}(z)^n=M_{\underline g}(z^n)$ 
for all $z\in S^1$ and 
$g\in \Diff(S^1)$.

$M^{(n)}$ is the unique isomorphism of $\Diff^{(n)}(S^1)$ into 
$\Diff(S^1)$ such that $M^{(n)}\restriction \text{\Mob}^{(n)}$ is 
given in Prop. 
\ref{cover}.
\end{proposition}
\proof
The proof is analogous to the proof of Prop. \ref{cover}.
\endproof

The Virasoro algebra is the infinite dimensional Lie algebra 
generated 
by elements $\{L_n \mid n\in\mathbb Z\}$ and $c$ with
relations
\begin{equation}\label{vir-rel}
[L_m,L_n]=(m-n) L_{m+n} + \frac{c}{12}(m^3-m)\delta_{m,-n}.
\end{equation}
and $[L_n,c]=0$. It is the (complexification of) the
unique, non-trivial one-dimensional central extension of the Lie 
algebra 
of $\Vect(S^1)$.

The elements $L_{-1},L_{0},L_{1}$ of the Virasoro 
algebra are clearly a basis of $s\ell(2,\mathbb C)$.
The Virasoro algebra  contains infinitely many
further copies of $s\ell(2,\mathbb C)$, generated by the elements
$L'_{-n},L'_{0},L'_{n},\; n>1$, where
\begin{gather}
L'_{n} \equiv \frac{1}{|n|}L_{n},\quad n \neq 0, \\
L'_{0} \equiv \frac{1}{n}L_{0} +  \frac{c}{24}\frac{(n^{2}-1)}{n}\ .
\label{spin}
\end{gather}         
For any fixed integer $n>0$ we have
\begin{equation}
  [L'_{n},L'_{-n}]=2L'_{0},\qquad
         [L'_{\pm n},L'_{0}]= \pm L'_{\pm n}
\end{equation}
which are indeed the relations for the usual generators in 
$s\ell(2,\mathbb C)$.

There is a one-to-one correspondence between projective irreducible 
unitary representations of $\Diff(S^1)$ and irreducible unitary 
representations of $\Diff^{(\infty)}(S^1)$.

We shall be interested in positive energy ($L_0 \geq 0$) 
representations of the Virasoro algebra which are 
unitary (i.e. preserving the involution). They correspond to 
projective unitary 
representations of $\Diff(S^1)$ with positive energy.

Given a projective unitary representation $U$ of $\Diff(S^1)$ and a 
fixed $n\in\mathbb N$, we 
obtain a projective unitary representation $U^{(n)}$ of 
$\Diff^{(n)}(S^1)$
\begin{equation}\label{ass}
U^{(n)}\equiv U\cdot M^{(n)}\ .
\end{equation}
(There is an analogous passage from unitary representations 
of $\Mob$ to unitary representations of $\Mob^{(n)}$.)

Starting with a positive energy, unitary representation $U$ of the 
Virasoro 
algebra with central charge $c$, it can be easily seen that the above 
construction (\ref{ass}) gives a positive energy, unitary 
representation $U^{(n)}$ of 
the Virasoro algebra with central charge $nc$. This will also be 
clear by the 
content of this paper. 

\section{Conformal nets on $S^1$}
\label{nets}

We denote by $\I$ the family of proper intervals of $S^1$. 
A {\it net} $\A$ of von Neumann algebras on $S^1$ is a map 
\[
I\in\I\to\A(I)\subset B(\H)
\]
from $\I$ to von Neumann algebras on a fixed Hilbert space $\H$
that satisfies:
\begin{itemize}
\item[{\bf A.}] {\it Isotony}. If $I_{1}\subset I_{2}$ belong to 
$\I$, then
\begin{equation*}
 \A(I_{1})\subset\A(I_{2}).
\end{equation*}
\end{itemize}
If $E\subset S^1$ is any region, we shall put 
$\A(E)\equiv\bigvee_{E\supset I\in\I}\A(I)$ with $\A(E)=\mathbb C$ 
if $E$ has empty interior (the symbol $\vee$ denotes the von Neumann 
algebra generated). 

The net $\A$ is called {\it local} if it satisfies:
\begin{itemize}
\item[{\bf B.}] {\it Locality}. If $I_{1},I_{2}\in\I$ and $I_1\cap 
I_2=\emptyset$ then 
\begin{equation*}
 [\A(I_{1}),\A(I_{2})]=\{0\},
 \end{equation*}
where brackets denote the commutator.
\end{itemize}
The net $\A$ is called {\it M\"{o}bius covariant} if in addition 
satisfies
the following properties {\bf C,D,E,F}:
\begin{itemize}
\item[{\bf C.}] {\it M\"{o}bius covariance}. 
There exists a strongly 
continuous unitary representation $U$ of $\Mob$ on $\H$ such that
\begin{equation*}
 U(g)\A(I) U(g)^*\ =\ \A(gI),\quad g\in\Mob,\ I\in\I.
\end{equation*}
\item[{\bf D.}] {\it Positivity of the energy}. The generator of the 
one-parameter rotation subgroup of $U$ (conformal Hamiltonian) is 
positive. 
\item[{\bf E.}] {\it Existence of the vacuum}. There exists a unit 
$U$-invariant vector $\Omega\in\H$ (vacuum vector), and $\Omega$ 
is cyclic for the von Neumann algebra $\bigvee_{I\in\I}\A(I)$.
\end{itemize}

The above axioms imply Haag duality (see \cite{BGL}): 
\[
\A(I)'=\A(I'),\quad I\in\I\ ,
\]
where $I'$ is the interior of $S^1\setminus I$.

\begin{itemize}
\item[{\bf F.}] {\it Irreducibility}. $\bigvee_{I\in\I}\A(I)=B(\H)$. 
Indeed $\A$ is irreducible iff
$\Om$ is the unique $U$-invariant vector (up to scalar multiples), 
and 
iff the local von Neumann 
algebras $\A(I)$ are factors. In this case they are III$_1$-factors 
(unless $\A(I)=\mathbb C$ identically), see \cite{GL2}.
\end{itemize}

By a {\it conformal net} (or diffeomorphism covariant net)  
$\A$ we shall mean a M\"{o}bius covariant net such that the following 
holds:
\begin{itemize}
\item[{\bf G.}] {\it Conformal covariance}. There exists a projective 
unitary representation $U$ of $\Diff(S^1)$ on $\H$ extending the unitary 
representation of $\PSL$ such that for all $I\in\I$ we have
\begin{gather*}
 U(g)\A(I) U(g)^*\ =\ \A(gI),\quad  g\in\Diff(S^1), \\
 U(g)xU(g)^*\ =\ x,\quad x\in\A(I),\ g\in\Diff(I'),
\end{gather*}
\end{itemize}
where $\Diff(I)$ denotes the subgroup of 
smooth diffeomorphisms $g$ of $S^1$ such that $g(z)=z$ for all $z\in 
I'$.

A representation $\pi$ of $\A$ on a Hilbert space $\H$ is a map 
$I\in\I\mapsto\pi_I$ that associates to each $I$ a normal 
representation of $\A(I)$ on $B(\H)$ such that
\[
\pi_{\tilde I}\res\A(I)=\pi_I,\quad I\subset\tilde I, \quad 
I,\tilde I\subset\I\ .
\]
$\pi$ is said to be M\"obius (resp. diffeomorphism) covariant if 
there is a projective unitary representation $U_{\pi}$ of $\Mob$ (resp. 
$\Diff^{(\infty)}(S^1)$) on $\H$ such that
\[
\pi_{gI}(U(g)xU(g)^*) =U_{\pi}(g)\pi_{I}(x)U_{\pi}(g)^*
\]
for all $I\in\I$, $x\in\A(I)$ and $g\in\Mob$ (resp. 
$g\in\Diff^{(\infty)}(S^1)$). Note that if $\pi$ is irreducible and 
diffeomorphism covariant then $U$ is indeed a projective unitary 
representation of $\Diff(S^1)$.

Following \cite{DHR}, given an interval $I$ and a representation $\pi$ 
of $\A$, there is an endomorphism of $\A$ localized in $I$ equivalent 
to $\pi$; namely $\r$ is a representation of $\A$ on the vacuum Hilbert 
space $\H$, unitarily equivalent 
to $\pi$, such that $\r_{I'}=\text{id}\restriction\A(I')$. We refer 
to \cite{GL2} for basic facts on this structure, in particular for 
the definition of the dimension $d(\r)$, that turns out to equal the 
square root of the Jones index \cite{L1}. The reader will also find 
basic notions concerning sectors of factors at the beginning of 
Sect. \ref{fus} or in \cite{Ko}.
\subsubsection{Restriction to the real line}
Denote by $\I_0$ the set of open, 
connected, non-empty, proper subsets of $\mathbb R$, thus $I\in\I_0$ 
iff $I$ is an open interval or half-line (by an interval of $\mathbb 
R$ we shall always mean a non-empty open bounded interval of $\mathbb 
R$).

Given a net $\A$ on $S^1$ we shall denote by $\A_0$ its restriction 
to $\mathbb R = S^1\setminus\{-1\}$. Thus $\A_0$ is an isotone map on
$\I_0$, that we call a \emph{net on $\mathbb R$}. 

A representation $\pi$ of $\A_0$ on a Hilbert space $\H$ is a map 
$I\in\I_0\mapsto\pi_I$ that associates to each $I\in\I_0$ a normal 
representation of $\A(I)$ on $B(\H)$ such that
\[
\pi_{\tilde I}\res\A(I)=\pi_I,\quad I\subset\tilde I, \quad 
I,\tilde I\in\I_0\ .
\]
A representation $\pi$ of $\A_0$ is also called a 
\emph{soliton}\footnote{There are more general soliton 
sectors, namely representations normal on left (resp. right) 
half-lines, but non-normal on right (resp. left) half-lines. These 
will not be considered in this paper.}.

Clearly a representation $\pi$ of $\A$ restricts to a soliton 
$\pi_0$ of $\A_0$. But a representation $\pi_0$ of $\A_0$ does not 
necessarily extend to a representation of $\A$. 

\subsection{Normality for $\a$-induction}
\label{3.1}
Let $\A$ be a M\"obius covariant net and $\B$ a subnet. Given  a 
bounded interval $I_0\in\I_0$ we fix canonical endomorphism 
$\g_{I_0}$ associated with $\B(I_0)\subset\A(I_0)$. Then we can choose 
for each $I\subset\I_0$ with $I\supset I_0$ a canonical endomorphism 
$\g_{I}$ of $\A(I)$ into $\B(I)$ in such a way that 
$\g_{I}\res\A(I_0)=\g_{I_0}$ and 
$\l_{I_1}$ is the identity on $\B(I_1)$ if 
$I_1\in\I_0$ is disjoint from $I_0$, where 
$\l_{I}\equiv\g_{I}\res\B(I)$.

We then have an endomorphism $\g$ of the $C^*$-algebra 
$\mathfrak A\equiv\overline{\cup_{I}\A(I)}$
($I$ bounded interval of $\mathbb R$).

Given a DHR endomorphism $\r$ of $\B$ localized in $I_0$, the 
$\a$-induction $\a_{\r}$ 
of $\r$ is the endomorphism of $\mathfrak A$ given by
\[
\a_{\r}\equiv \g^{-1}\cdot\Ad\e(\r,\l)\cdot\r\cdot\g\ ,
\]
where $\e$ denotes the right braiding unitary symmetry (there is 
another choice for $\a$ associated with the left braiding).
$\a_{\r}$ is localized in a right half-line containing $I_0$, namely 
$\a_\r$ is the identity on $\A(I)$ if $I$ is a bounded interval 
contained in the left complement of $I_0$ in $\mathbb R$. Up to 
unitarily equivalence, $\a_\r$ is localizable in any right half-line 
thus $\a_\r$ is normal on left half-lines, that is to say, for every 
$a\in\mathbb R$, 
$\a_\r$ is normal on the $C^*$-algebra 
$\mathfrak A(-\infty,a)\equiv\overline{\cup_{I\subset 
(-\infty,a)}\A(I)}$
($I$ bounded interval of $\mathbb R$), namely $\a_\r\res\mathfrak 
A(-\infty,a)$ extends to a normal morphism of $\A(-\infty,a)$. 

We now show that $\a_\r$ is normal on right half-lines. To this end we 
use the fact that our nets on $\mathbb R$ are restrictions of nets 
on $S^1$.
\begin{proposition}
$\a_\r$ is a soliton endomorphism of $\A_0$.
\end{proposition}
\proof
It is convenient to use the circle picture, thus $I_0\subset 
S^1\setminus\{-1\}$, say $I_0=(a,b)$ where $a,b\in 
S^1\setminus\{-1\}$, and $b>a$ in the counterclockwise order 
(intervals 
do not contain $-1$). Let $a_n,b_n\in S^1\setminus\{-1\}$ with $a < 
b < a_n < b_n$ 
and $\r_n$ an endomorphism of $\B$ equivalent to $\r$ and localized 
in $(a_n ,b_n)$. 
With $u_n\in\B(a,b_n)$ a unitary such that $\r_n=\Ad u_n\cdot\r$, we 
have
\[
\a_\r\res \A(c, a_n) = \Ad u^*_n \res \A(c, a_n),
\]
for every $c<a$. Going to the limit $c\to -1^-$, $b_n\to -1^+$ the 
above gives the definition of $\a_\r$ on the $C^*$-algebra $\mathfrak 
A$ originally given in \cite{R}.

We want to show that $\a_\r\res\mathfrak {\mathfrak A}(d,-1)$ extends 
to a normal map of $\A(-1,d)$ for any given $d\neq -1$.

Now, as $\B$ is defined on $S^1$, we may push the interval $(a_n,b_n)$ 
even beyond the point $-1$. Namely we may choose an interval $(a',b')$ 
with $-1<a'<b'<a$, an 
endomorphism $\r'$ of $\B$ equivalent to $\r$ localized in $(a',b')$, 
and a unitary $u\in\B(a,b')$ such that $\r'=\Ad u\cdot\r$. 
Then $\a_\r\res \mathfrak A(a, -1) =
\Ad u^* \res \mathfrak A(a, -1)$, showing that $\a_\r$ extends to a 
normal morphism of $\A(a, -1)$. Of course we 
may take a smaller $a$ in the definition of $I_0$, thus 
$\a_\r$ is normal on all right half-lines.
\endproof
\subsection{CMS property}
In this section $\A$ is a M\"obius covariant local net on $S^1$.
We shall say that $\A$ has \emph{property CMS} if it admits at most 
countably many different irreducible (DHR) sectors and all of them 
have finite index.

Let $\b$ be a vacuum preserving, involutive automorphism of $\A$ 
and $\B=\A^{\b}\subset\A$ the fixed-point 
subnet. The restriction of $\A$ and $\B$ to $\mathbb R = 
S^1\setminus\{-1\}$ are denoted by $\A_0$ and $\B_0$ as above.

We denote by $[\s]$ the sector of $\B$ dual to $\b$. Choosing 
an interval $I_0\subset\mathbb R$ there is a unitary 
\[
v\in\A(I),\: v^*=v,\: \b(v)= -v \ .
\]
Then $\s\equiv\Ad v\res \B$ is an automorphism of $\B$ localized in 
$I_0$. We have $d(\s)=1$ and $\s^2 = 1$.

Given a DHR endomorphism $\mu$ of $\B$ localized in an interval 
$I_0\subset\mathbb R$, we denote as above by $\a_{\mu}$ the 
right $\a$-induction of $\mu$ to $\A_0$. Recall that in general $\a_{\mu}$ is a 
soliton sector of $\A_0$. With $\e(\mu,\s)$ the right 
statistics operator, the condition for $\a_{\mu}$ to be localized 
in a bounded interval of $\mathbb R$, i.e. to be a DHR endomorphism 
of $\A$, is that the monodromy operator 
$\e(\mu,\s)\e(\s,\mu)=1$. If $\mu$ is localized left to 
$\s$, then $\e(\s,\mu)=1$, so we have:
\begin{proposition}
Let $\mu$ be localized in an interval $I\subset \mathbb R$ in the 
left 
complement of $I_0$ in $\mathbb R$. Then
$\a_{\mu}$ is a DHR sector of $\A$ iff $\e(\mu,\s)=1$.
\end{proposition}
\proof
Let $I_1$ be an interval of $\mathbb R$ in the right complement of 
$I_0$, $\mu'$ an endomorphism of $\B$ localized in $I_1$ and 
$u\in\Hom(\mu,\mu')$ a unitary. Then $\mu(x)=\Ad u^*(x)$ for all 
$x\in\B(I_2)$ if $I_2$ is an interval left to $I_1$. 

We then have $\a_{\mu}(x)=\Ad u^*(x)$ if $x\in\A(I_2)$. It follows 
that $\a_{\mu}$ is localized in $I_1$ iff $\a_{\mu}$ acts trivially on
$\A(I_0)$. As $\A(I_0)$ is generated by $\B(I_0)$ and $v$, this is 
the 
case iff 
\[
\a_{\mu}(v)=v \Leftrightarrow u^*vu = v  \Leftrightarrow 
\varepsilon(\mu,\s)= u^*\s(u) = 1\ .
\]
\endproof
Let $\mu$ be an irreducible endomorphism localized left to $\s$. 
As $\e(\mu,\s)\in\Hom(\mu\s,\s\mu)$ and $\s$ and $\mu$ commute, it 
follows that $\e(\mu,\s)$ is scalar. Denoting by $\iota$ the identity 
sector, by the braiding fusion relation we have
\[
1=\e(\mu,\iota)=\e(\mu,\s^2)
=\s(\e(\mu,\s))\e(\mu,\s)=\e(\mu,\s)\e(\mu,\s)\ ,
\]
thus $\e(\mu,\s)=\pm 1$.

If $\mu$ is not necessarily irreducible, we shall say that $\mu$ is 
\emph{$\s$-Bose} 
if $\e(\mu,\s)= 1$ and that $\mu$ is \emph{$\s$-Fermi} 
if $\e(\mu,\s)= -1$. As we have seen, if $\mu$ is irreducible then 
$\mu$ is either $\s$-Bose or $\s$-Fermi.
\begin{corollary}\label{munu}
Let $\mu,\nu$ be DHR sectors of $\B$. If $\mu,\nu$ are 
both $\s$-Fermi, then $\a_{\mu\nu}$ is a DHR sector of $\A$.
\end{corollary}
\proof
We may assume that both $\mu$ and $\nu$ are localized left to $\s$.
By the braiding fusion relation we have 
\[
\e(\mu\nu,\s)=\mu(\e(\nu,\s))\e(\mu,\s)=\e(\nu,\s)\e(\mu,\s)=1\ .
\]
\endproof
\begin{lemma}\label{dirint}
Let $\mu$ be a $\s$-Bose sector of $\B$. Then $\mu$ has a direct 
integral decomposition into irreducible $\s$-Bose sectors.
\end{lemma}
\proof
$\a_{\mu}$ is a $\s$-Bose sector of $\A$, thus $\a_{\mu}$ has a 
direct integral decomposition into irrreducibles \cite{KLM}, say 
$\a_{\mu}=\int^{\oplus}\pi_{t}\text{d}m(t)$. Since $\B\subset\A$ is 
a finite-index subnet, the restriction of $\pi_t$ 
to $\B$ is the sum of finitely many irreducible $\s$-Bose 
representations, so the restriction of $\a_{\mu}$ to $\B$ has a direct
integral decomposition into irreducible $\s$-Bose sectors.  By
Frobenius reciprocity (cf. Th. B.2) $\mu$ is contained in the restriction of
$\a_{\mu}$ to $\B$ and we are done.
\begin{corollary}\label{AB-CMS}
Assume $\A$ to have property CMS, then $\B$ has property CMS.
\end{corollary}
\proof
First suppose that $\B$ has an irreducible $\s$-Bose sector $\mu$ 
with $d(\mu)=\infty$.
Then $\a_{\mu}$ is a DHR sector of $\A$ with 
$d(\a_{\mu})=d(\mu)=\infty$. 

As $\A$ has property CMS, there is an irreducible finite-index DHR 
sector $\l$ of $\A$ with $\l\prec\alpha_{\mu}$.

By Frobenius reciprocity we have the equality between the the 
dimensions of the intertwiners spaces 
$\langle \a_{\mu}, \l\rangle = \langle\mu, \g\l\res\B\rangle$, thus 
$\g\l\res\B\succ\mu$. As $d(\g\l\res\B)<\infty$ 
then $d(\mu)<\infty$ and this shows that $\B$ has no irreducible 
$\s$-Bose sector with infinite dimension.

Suppose now that $\B$ has uncountably many $\s$-Bose irreducible 
sectors 
$\{\mu_i\}$ with finite dimension. As $\A$ has property CMS there 
must be an irreducible finite dimensional DHR sector $\l$ such that 
$\a_{\mu_i}\succ\l$ for uncountably many $i$. By Frobenius reciprocity
$\mu_i\prec \g\l\res\B$, thus 
$d(\g\l\res\B)=\infty$, which is not possible because 
$d(\l)<\infty$. Thus $\A$ admits at most countably many inequivalent 
irreducible $\s$-Bose sectors and all have finite dimension.

Suppose now that $\mu$ an is irreducible, $\s$-Fermi and infinite 
dimensional 
sector of $\B$. Then $\bar\mu\mu$ is $\s$-Bose. Now $\B$ inherits the 
split property from $\D$ (this is rather immediate, see \cite{L4}) so 
$\bar\mu\mu$ has a direct integral decomposition into irreducible 
sectors that must be almost everywhere $\s$-Bose because $\bar\mu\mu$ 
is $\s$-Bose.

By what we have proved above, 
$\bar\mu\mu$ is then a direct sum of finite dimensional $\s$-Bose 
sectors, and analogously the same is  true for $\mu\bar\mu$, and this 
entails $d(\mu)<\infty$ as in of Lemma \ref{finindex}.

It remains to show that $\B$ cannot have uncountably many $\s$-Fermi 
irreducible sectors 
$\{\mu_i\}$ with finite dimension. On the contrary for a given $i_0$ 
there should exist uncountably many $i$ and a fixed finite dimensional 
irreducible sector $\l$ of $\B$ such that $\mu_{i_0}\mu_i\succ \l$ 
because we have already proved that there are at most countably many 
finite dimensional $\s$-Bose irreducible sectors. 
By Frobenius reciprocity then $\bar\l\mu_{i_0}\succ\bar\mu_i$, which is 
not possible because $d(\mu_{i_0})<\infty$.

This concludes our proof.
\endproof

\begin{lemma}\label{finindex}
Let $M$ be a factor and $\r\in\End(M)$ an irreducible endomorphism. 
If 
there are $\s,\s'\in\End(M)$ such that $\r\s\succ\mu$ and 
$\s'\r\succ\mu'$ with $\mu,\mu'$ 
finite index subsectors, then $d(\r)<\infty$.
\end{lemma}
\proof
With $\r'\equiv \s\bar\s\bar\r\in\End(M)$, we have 
$\r\r'=\r\s\bar\s\bar\r\succ\mu\bar\mu\succ\iota$. Analogously there 
is 
$\r''\in\End(M)$ such that $\r''\r\succ\iota$, thus $\r$ has finite 
index by the criterion on the existence of the conjugate sector in \cite{L1'}.
\endproof
\section{Canonical representation of $\A_0\otimes\A_0$}
For simplicity we shall now consider the $2$-fold tensor product, 
which is however sufficient for most of the applications.
We shall return on this point in later sections and have a more 
general analysis in the case of arbitrary $n$-fold tensor product. 

We shall say that a set $E\subset S^1$ is a \emph{symmetric 
2-interval} if $E=I_1\cup I_2$ where $I_1,I_2\in\I$ are interval of 
with 
length less than $\pi$ and $I_2=R(\pi)I_1=-I_1$. The set of all 
symmetric 2-intervals is denoted by $\I^{(2)}$. 

Given an interval $I\in\I$, then 
\[
E\equiv\sqrt{I}=\{z\in S^1\, |\, z^2\in I\}
\]
is a symmetric 2-interval. Conversely, given a symmetric 2-interval 
$E=I_1\cup I_2 $, 
then $I\equiv I_1^2 =I_2^2$ is an interval and $E=\sqrt{I}$, 
thus there is a bijection between $\I$ and $\I^{(2)}$.

In the following $\A$ denotes a diffeomorphism covariant, local net 
of von Neumann algebras on $S^1$. We denote by $U$ the associated projective 
unitary representation of $\Diff(S^1)$. 

We assume the \emph{split property}. 

Given $\zeta\in S^1$, we shall denote by $\sqrt[(\zeta)]{\cdot}$ 
the square root function on $S^1$ with a discontinuity in $\zeta$, 
namely $z\in S^1\mapsto \sqrt[(\zeta)]{z}\in S^1$ is the unique 
function such that $(\sqrt[(\zeta)]{z})^2 =z$, $\sqrt[(\zeta)]{1}=1$, 
$\sqrt[(\zeta)]{\cdot}$ is continuous at all $z\neq\zeta$ and 
continuous from the right (counterclockwise) at $z=\zeta$.

Let $I\subset S^1$ be an interval and set $E=\sqrt{I}\in\I^{(2)}$, 
$E=I_1\cup I_2$.
Given $\zeta\in \overline{I'}$ 
we choose the two components of $E$ so that $I_1 = 
\sqrt[(\zeta)]{I}$, 
$I_2 = -\sqrt[(\zeta)]{I}=R(\pi)I_1$.

Let $h\in\Diff(S^1)$ 
be such that $h(z)=\sqrt[(\zeta)]{z}$, $z\in I$ (cf. \cite{S}) and 
set $\bar h(z)= -h(z)$. Clearly $\bar h\in\Diff(S^1)$ and 
$\bar h(z)=-\sqrt[(\zeta)]{z}$, $z\in I$.
Setting
\begin{gather}
\Phi^{(\zeta)}_{I} \equiv \Ad U(h)\res \A(I), \\
\bar\Phi^{(\zeta)}_{I} \equiv \Ad U(\bar h)\res \A(I),
\end{gather}
by diffeomorphism covariance $\Phi^{(\zeta)}_{I}, 
\bar\Phi^{(\zeta)}_{I}$
are isomorphisms of $\A(I)$ with $\A(I_1)$ and with $\A(I_2)$.
\begin{proposition}\label{zz'}
Let $I\in\I$ and $\z\notin I$. We have:
\begin{itemize}
\item[$(a)$] $\Phi^{(\zeta)}_{I},\bar\Phi^{(\zeta)}_{I}$ do not 
depend on the choice of $h$.
\item[$(b)$] $\bar\Phi^{(\zeta)}_{I}=\Ad U(R(\pi))\cdot 
\Phi^{(\zeta)}_{I}$.
\item[$(c)$] If $\z'\notin I$, then 
$\Phi^{(\zeta')}_{I}=\Phi^{(\zeta)}_{I}$ or 
$\Phi^{(\zeta')}_{I}=\bar\Phi^{(\zeta)}_{I}$.
Denote by $[\z,\z')$ the interval of $S^1$ in the counterclockwise
order and assume $(\z,\z')\nsubseteqq I$. Then 
$\Phi^{(\zeta)}_{I}=\Phi^{(\zeta')}_{I}$ iff $1\notin (\z,\z')$.
\end{itemize}
\end{proposition}
\proof
$(a)$: Let $k\in\Diff(S^1)$ be such that $k\res I 
=\sqrt[(\z)]{\cdot}$. Then $ k^{-1}\cdot h\res I$ is the 
identity, thus $V\equiv U(k^{-1}\cdot h)\in\A(I')$ and $\Ad 
U(h)\res\A(I)=\Ad U(k)V\res\A(I)=\Ad U(k)\res\A(I)$.

$(b)$: We have
\begin{multline}
\bar\Phi^{(\zeta)}_{I}=\Ad U(\bar h)\res \A(I)=
\Ad U(R(\pi)\cdot h)\res \A(I)\\
=\Ad U(R(\pi))\Ad U(h)\res \A(I)
=\Ad U(R(\pi))\cdot \Phi^{(\zeta)}_{I}.
\end{multline}

$(c)$: The restriction of $h$ to $I$ does not vary as long we 
choose another $\z'\notin I$ such that
$\sqrt[(\z')]{z}=\sqrt[(\z)]{z}$ for all $z\in I$, thus  
$\Phi^{(\zeta')}_{I}=\Phi^{(\zeta)}_{I}$ for such $\z'$. Otherwise 
$h(z)$ changes to $-h(z)$, $z\in I$, and then 
$\Phi^{(\zeta')}_{I}=\bar\Phi^{(\zeta)}_{I}$. The rest is now clear. 
\endproof
We now set 
$\pi^{(\zeta)}_{I}\equiv 
\chi_{I}\cdot(\Phi^{(\zeta)}_{I}\otimes\bar\Phi^{(\zeta)}_{I})$,
where $\chi_I$ is the canonical isomorphism of 
$\A(I_1)\otimes\A(I_2)$ 
with $\A(I_1)\vee\A(I_2)$ given by the split property. In other words 
$\pi^{(\zeta)}_{I}$ is the unique isomorphism of $\A(I)\otimes\A(I)$ 
with $\A(I_1)\vee\A(I_2)$ such that
\begin{equation}\label{piz}
\pi^{(\zeta)}_{I}(x_1\otimes x_2)=
\Phi^{(\zeta)}_{I}(x_1)\bar\Phi^{(\zeta)}_{I}(x_2), \quad x_1 , x_2 
\in \A(I)\ .
\end{equation}
\begin{proposition}\label{ext}
Let $I\subset\tilde I$ be intervals and $\zeta\notin\tilde I$. Then 
$\pi^{(\zeta)}_{\tilde I}\res 
\A(I)\otimes\A(I)=\pi^{(\zeta)}_{I}$.
\end{proposition}
\proof
Immediate by the above Proposition \ref{zz'}.
\endproof
\begin{corollary}\label{zz'2}
Let $I$ be an interval and $\zeta,\zeta'\notin I$.  Then either 
$\pi^{(\zeta')}_{I}=\pi^{(\zeta)}_{I}$ or 
$\pi^{(\zeta')}_{I}=\pi^{(\zeta)}_{I}\cdot\a$, where $\a$ is the flip 
automorphism of $\A\otimes\A$.

The first alternative holds iff $\zeta,\zeta'$ both belong or both do 
not belong to the closure of the connected component of 
$I'\setminus\{1\}$ intersecting the upper half plane.
\end{corollary}
\proof
If $\Phi^{(\zeta')}_{I}=\Phi^{(\zeta)}_{I}$, then also 
$\bar\Phi^{(\zeta')}_{I}=\bar\Phi^{(\zeta)}_{I}$ and then clearly 
$\pi^{(\zeta')}_{I}=\pi^{(\zeta)}_{I}$. 

Otherwise $\Phi^{(\zeta')}_{I}=\bar\Phi^{(\zeta)}_{I}$, thus 
$\Phi^{(\zeta)}_{I}=\bar\Phi^{(\zeta')}_{I}$, and we have 
\begin{multline}
\pi^{(\zeta)}_{I}(\a(x_1\otimes x_2))=
\Phi^{(\zeta)}_{I}(x_2)\bar\Phi^{(\zeta)}_{I}(x_1)
=\Phi^{(\zeta')}_{I}(x_1)\bar\Phi^{(\zeta')}_{I}(x_2)\\
=\pi^{(\zeta')}_{I}(x_1\otimes x_2),
\quad x_1 , x_2 \in \A(I)\ .
\end{multline}
The rest follows by Prop. \ref{zz'}.
\endproof
In the following we shall denote the net $\A\otimes\A$ by $\D$.

As usual we may identify $S^1\setminus\{-1\}$ with $\mathbb R$ by the 
stereographic map. Let $\A_0$ be the net on $\mathbb 
R$ obtained by restricting $\A$ to $S^1\setminus\{-1\}$. We denote 
by $\pi$ the restriction of  $\pi^{(\zeta=-1)}$ to $\D_0=\A_0\otimes\A_0$.
\begin{proposition}\label{irr0}
$\pi$ is a representation of $\D_0$. Indeed $\pi$ is 
an irreducible soliton.
\end{proposition}
\proof
That $\pi$ is a soliton representation follows from Prop. \ref{ext}
and the fact that $\pi^{(\z)}$ is normal on $\D(I)$ for every interval $I$ 
not containing $\z$, including the case $\z\in\bar I$ (half-lines).

Now $\pi(\D_0(I))=\A(E)$ where $E=\sqrt{I}$, thus 
\[
\bigvee_{\z\notin I\in\I}\pi(\D_0(I)) = \bigvee_{\pm i\notin 
E\in\I^{(2)}}\!\!\!\A(E) = \A(S^1\setminus\{i,-i\}) = B(\H)
\]
because $\A$ is 2-regular by Haag duality and the factoriality of the 
local von Neumann algebras, so $\pi$ is irreducible.
\endproof
By Prop. \ref{irr0} $\pi$ is a representation of $\D_0$, namely 
$\pi$ is consistently defined on all von Neumann algebras 
$\D_0$ with $I\subset\mathbb R$ either an interval or an 
half-line. However $\pi$ is \emph{not} a DHR representation of $\D_0$ 
namely, given an interval $I_0\subset \mathbb R$, $\pi$ is not normal on the 
$C^*$-algebra $\mathfrak D(I'_0)\equiv\overline{\cup_{I\subset I'_0}\D(I)}$ 
($I$ interval of $\mathbb R$).

As $\D_0$ satisfies half-line duality, namely
\[
\D_0 (-\infty,a)'=\D_0 (a,\infty), \quad a\in\mathbb R,
\]
by the usual DHR argument \cite{DHR} $\pi$ is unitarily equivalent to 
a representation $\r$ of $\D_0$ on $\cal H\otimes\cal H$ which 
acts identically on $\D_0(-\infty,0)$, thus $\r$ restricts 
to an endomorphism of $\D_0(0,\infty)$.
\begin{proposition}\label{iso2}
Setting $M\equiv\D_0(0,\infty)$, the 
inclusion $\r(M)\subset M$ is isomorphic to the 2-interval 
inclusion $\A(E)\subset\hat\A(E)$.
\end{proposition}
\proof
In the circle picture with $\z= -1$, setting $I=S^+$ (the upper semicircle)
and $E\equiv \sqrt{I}$, we have $M=\D(I)$ and
\[
\A(E)=\pi_I^{(\z)}(\D(I)),\quad 
\A(E')=\pi_{I'}^{(\z)}(\D(I')),
\]
thus we have an equality of inclusions:
\[
\big \{\A(E)\subset\hat\A(E)\big \} =
\big \{\pi_I^{(\z)}(\D(I))\subset\pi_{I'}^{(\z)}(\D(I'))'\big \}\ .
\]
As $\pi$ is unitarily equivalent to $\r$ and $\r_{I'}$ is the 
identity 
on $\D(I')$, the second inclusion is isomorphic to 
\begin{equation}
\big \{\r_I(\D(I))\subset\r_{I'}(\D(I'))'\big \}
=
\big \{\r_I(\D(I))\subset \D(I')'\big \}
=
\big \{\r_I(\D(I))\subset\D(I)\big \}\ .
\end{equation}
\endproof

\subsection{Canonical representation of $(\A\otimes\A)^{\text{flip}}$}

We shall denote by $\B\equiv(\A\otimes\A)^{\a}$ the fixed-point 
subnet of $\D$ with 
respect to the flip symmetry $\a$.

By Prop. \ref{zz'2} 
$\pi^{(\zeta)}_{I}\res\B(I)=\pi^{(\zeta')}_{I}\res\B(I)$ 
for all $\z,\z'\notin I$, therefore
\[
\t_{I}\equiv \pi^{(\zeta)}_{I}\res\B(I)
\]
is independent of $\z\notin I$ and thus well defined.

Recall that the spin of a M\"{o}bius covariant representation is the 
lowest eigenvalue of the conformal Hamiltonian $L_0$ in the 
representation space.
\begin{corollary}\label{tirr}
$\t : I\mapsto \t_{I}$ is a (DHR) diffeomorphism covariant 
representation of $\B$ (with positive energy).
The covariance unitary representation is given by 
\[
U^{(2)}\equiv U\cdot M^{(2)}
\]
(see Sect. \ref{first}), where $U$ is the covariance 
unitary representation associated with $\A$.

$\t$ is direct sum of two irreducible diffeomorphism covariant 
representations with spin $c/16$  and $1/2 + c/16$.
\end{corollary}
\proof
It follows by Prop. \ref{irr0} that $\t$ is a representation.

We shall show that the projective unitary representation 
$U^{(2)}\equiv U\cdot 
M^{(2)}$ of $\Diff^{(2)}(S^1)$ implements the covariance of $\t$, 
namely, setting $\bar U(g)\equiv U(g)\otimes U(g)$,
\[
\t_{gI}(\bar U (g)x\bar U(g)^*) = U^{(2)}(g)\t_{I}(x)U^{(2)}(g)^*, \: I\in\I, 
x\in\B(I), g\in\Diff^{(2)}(S^1) \ .
\]
The above formula will follow if we show that 
\[
\pi^{(\z')}_{gI}(\bar U (g)x\bar U (g)^*) = 
U^{(2)}(g)\pi^{(\z)}_{I}(x)U^{(2)}(g)^*, \: I\in\I, 
x\in\D(I), g\in\Diff^{(2)}(S^1),
\]
for some $\z\notin I,\z'\notin gI$, and indeed it will suffice to verify 
this for $x=x_1\otimes 1$ or $x=1\otimes x_2$, $x_1 ,x_2\in\A(I)$. Suppose 
$x=x_1\otimes 1$:
\begin{multline}
U^{(2)}(g)\pi^{(\z)}_{I}(x)U^{(2)}(g)^*=\Ad 
U^{(2)}(g)\Phi^{(\z)}_{I}(x_1)
= \Ad U^{(2)}(g) U(h)(x_1)\\
=\Ad  U(h_g) 
U(\underline{g})(x_1)=\pi^{(g\z)}_{gI}(U(g)x_1 U(g)^*)=
\pi^{(g\z)}_{gI}(\bar U (g)x\bar U (g)^*)
\end{multline}
where $h(z)=\sqrt{z}$ on $I$ and $h_g(z)=\sqrt{z}$ on $gI$ (see Prop. 
\ref{cover2}). 

The computation in the case $x=1\otimes x_2$ is analogous.

Concerning the last statement, set $\mathfrak B_0(\mathbb R)$ for
the $C^*$-algebra $\overline{\cup_I \B(I)}$ ($I$ bounded interval) 
and analogously for $\mathfrak D_0$ and note that
\[
\mathbb C = \pi(\mathfrak D_0(\mathbb R))'= \{\t(\mathfrak 
B_0(\mathbb R)), 
\pi(v)\}' \ .
\]
Thus $\Ad\pi(v)$ acts ergodically on $\t(\mathfrak B_0(\mathbb 
R))'$. 
Since $v^2 =1$, $\text{dim}(\t(\mathfrak B_0(\mathbb R))')\leq 2$.
As $U^{(2)}(R(2\pi))=U(R(\pi))$ belongs to $\t(\mathfrak B_0(\mathbb 
R))'$, we have 
$\text{dim}(\t(\mathfrak B_0(\mathbb R))')= 2$, thus $\t$ has exactly 
two irreducible direct summands.

The spin of these two representations is now soon computed by formula
(\ref{spin}).
\endproof
\begin{lemma}\label{CMS2}
Let $\A$ be a local, split conformal net with the CMS property. Then  
$\mu_{\A}<\infty$.
\end{lemma}
\proof
The CMS property holds for $\D$ by \cite{KLM} (irreducible sectors 
of $\D$ are tensor product of irreducible sectors of $\A$). Thus for
$\B$ by Cor. \ref{AB-CMS}.

Now $\t$ is the sum of two irreducible 
representations, thus, by the CMS property, $\t$ has finite index.

With $I=S^+$ and $E=\sqrt{I}$ we have:
\begin{gather}
\A(E)=\pi_I(\D(I))\supset \t_I(\B(I))\\
\A(E')=\pi_{I'}(\D(I'))\supset \t_{I'}(\B(I'))
\end{gather}
thus 
\[
\t_I(\B(I))\subset \A(E)\subset \hat\A(E)\subset \t_{I'}(\B(I'))'\ ,
\]
but $\t_I(\B(I))\subset\t_{I'}(\B(I'))'$ has finite index and this 
entails $[\hat\A(E):\A(E)]<\infty$.
\endproof
Recall that a local net $\A$ is said to be \emph{$n$-regular} if 
$\mathfrak 
A(S^1\setminus F)$ is irreducible if $F\subset S^1$ is a finite set 
with $n$ points, namely $(\vee_{I\cap F=\emptyset}\A(I))'=\mathbb C$ 
($I\in\I$). 

It is immediate that, if $\A$ is conformal, $n$-regularity does 
not depend on the choice of the $n$-point $F$ set and
\begin{equation}
\A\ \text{is $2n$-regular} \Leftrightarrow \A(E)\subset\hat\A(E)\ 
\text{is 
irreducible}
\end{equation} 
where $E$ is any $n$-interval.
\begin{lemma}\label{irr2intincl}
Let $\A$ be a local, split conformal net. If $\mu_\A < \infty$ then 
the 2-interval inclusion $\A(E)\subset\hat\A(E)$ is irreducible. 
Thus $\A$ is 4-regular.
\end{lemma}
\proof
By Prop. \ref{iso2} we have to show that 
$\r_I(\D(I))'\cap\D(I)=\mathbb C$. 
This would follow from the theorem on 
the equivalence between local and global 
intertwiners \cite{GL2}, but $\r$ is not a DHR representation and that 
theorem does not apply here directly, but it will nevertheless give 
the result.

Let $T\in\r_I(\D(I))'\cap\D(I)= \big 
(\r_I(\D(I))\vee\r_{I'}(\D(I'))\big )'$.
Then $T\in \th_I(\B(I))\vee\th_{I'}(\B(I'))$, thus $T\in 
\th(\mathfrak B)'$ due to the equivalence between local and global 
intertwiners, because $\r\res \B$ is a covariant, finite-index 
representation. 

On the other and $T\r_I(v)=\r_I(v)T$, thus $T$ commutes with 
$\{\r_{\tilde I}(\B(\tilde I)), \r_{\tilde I}(v)\}''
=\r_{\tilde I}(\D(\tilde I))$ for all intervals $\tilde 
I\supset I$ and $T$ is a scalar because $\r$ is irreducible.
\endproof
We now state and begin to prove the dichotomy.
\begin{theorem} 
Let $\A$ be a local conformal net with the split 
property. Assume that every irreducible sector of $\A$ is finite 
dimensional. We then have the following dichotomy: Either
\begin{itemize}
\item[$(a)$]  $\A$ is completely rational or 
\item[$(b)$] $\A$ has uncountably many different irreducible sectors.
\end{itemize}
\end{theorem}
\proof
Assuming that $\A$ has the CMS property, we have to show that $\A$ is 
completely rational. By Lemma \ref{CMS2} $\mu_{\A}<\infty$, thus we 
have to show that, for a local conformal net with the split property, 
the implication ``$\mu_{\A}<\infty\Rightarrow$ strong additivity'' 
holds. This will be the content of Section \ref{stradd}.
\endproof
\subsection{The canonical endomorphism of the $n$-interval inclusion}
\label{4.2}
In this section $\A$ is again a local conformal net with the split property and 
$\D=\A\otimes\A$. Our results have direct extension to the case of a 
general $n$-fold tensor product, but we deal with the case $n=2$ for 
simplicity, but in the last corollary.

We keep the above notations, thus
$\pi$ is the canonical representation of $\D_0$ and $\r$ is a 
soliton endomorphism of $\D_0$ equivalent to $\pi$ and
localized in $S^+$. The conjugate sector 
$\bar\r$ of $\r$ is given by $[\bar\r] = [j\cdot\r\cdot j]$ where $j=\Ad J$ 
with $J$ the modular conjugation of $(\A(S^+),\Om)$ \cite{GL2}. Note that 
$j\cdot\r\cdot j$ is localized in the lower semicircle
$S^-$ but, as $\bar\r$ is normal on 
$\A(S^-)$, we can choose, in the same unitary equivalence class of 
$j\cdot\r\cdot j$, an endomorphism $\bar\r$ localized in $\A(S^+)$.
\begin{proposition}\label{4.11}
$\bar\r\r$ is a soliton of $\D_0$ localized in $S^+$.
\end{proposition}
\proof
The statement is clear by the above comments, as both $\r$ and 
$\bar\r$ are solitons localized in $S^+$.
\endproof
Denote by $\l_E$ the dual canonical endomorphism associated with the 
inclusion $\A(E)\subset\hat\A(E)$.
\begin{proposition}\label{intercan}
Let $\r$ be localized in the right half-line   
$I\subset \mathbb R\simeq S^1\setminus\{-1\}$. 
If $S^1\setminus\{-1\}\supset \tilde I\supset I$ is a half-line, the two 
squares of inclusions
\begin{equation*}
\begin{matrix} \D(I) & \subset & \D(\tilde I)\\ 
\cup & {} & \cup \\ 
\r_I(\D(I)) & \subset & \r_{\tilde I}(\D(\tilde I))
\end{matrix}
\ \qquad \text{and}\ \qquad
\begin{matrix} 
\hat\A(E) & \subset & \hat\A(\tilde E)\\
\cup & {} & \cup \\
\A(E) & \subset & \A(\tilde E)
\end{matrix}
\end{equation*}
are isomorphic, where $E=\sqrt{I}$, $\tilde E=\sqrt{\tilde I}$.

If $\bar\r$ is also localized in $I$   
the isomorphism $\pi_{\tilde I}:\D(\tilde I)\to\A(\tilde E)$, 
interchanges $[\bar\r_I\r_I\res\D(I)]$ and $[\l_E]$.
\end{proposition}
\proof
Let $U$ be a unitary from $\H$ to $\H\otimes\H$ such that 
$\pi_{I'}=\Ad U\res\D(I')$. Then we can assume $\r_I = \Ad 
U^*\cdot\pi_I$. 
The isomorphism $\pi_I:\D(I)\to\A(E)$ is thus the composition
\begin{equation}\label{cd}
\CD
\D(I)@>\r_I>>
\r_I(\D(I))   @>\Ad U> 
>  \A(E)
\endCD
\end{equation}
$\Ad U$ maps $\r_I(\D(I))$ onto $\A(E)$ and $\D(I)$ onto 
$\hat\A(E)$ as in Prop. \ref{iso2}. As $\tilde I'\subset I'$, we also 
have $\pi_{\tilde I'}=\Ad 
U\res\D(\tilde I')$, therefore $\Ad U$ maps $\r_{\tilde I}(\D(\tilde 
I))$ 
onto $\A(\tilde E)$ and $\D(\tilde I)$ onto 
$\hat\A(\tilde E)$, thus $\Ad U$ implements an isomorphism between 
the two squares.

In particular $\Ad U$ will interchange $\l_E$ with the dual canonical 
endomorphism associated with $\r_I(\D(I))\subset\D(I)$, which is 
$\r_I\bar\r_I\res \r_I(\D(I))$ (here $\bar\r_I$ is the 
conjugate  of $\r_I$ as sectors of $\D(I)$). Then $\r_I$ will 
interchange the latter with $\r_I^{-1}\r_I\bar\r_I\r_I =\bar\r_I\r_I$.
\endproof
It will follow from the results in Sect. \ref{tts} that, in the 
case $n=2$, $\r$ is self-conjugate, 
as both $\r$ and $\bar\r$ are associated with a degree 2 map on $S^1$. 
In the case of the $n$-fold tensor product this fact is not any 
longer true and we shall have a formula for $\bar \r$ in Prop. \ref{pi} 
which gives 
\begin{equation}\label{conj}
\bar \r \simeq \b_p^{-1}\cdot\r\cdot\b_p\ ,
\end{equation}
where $\b$ is the natural action of $\mathbb P_n$ on 
$\A\otimes\cdots\otimes\A$ and $p\in\mathbb P_n$ is the inverse map on 
the group $\mathbb Z_n$.

As a corollary of Prop. \ref{intercan} 
we now show that in the completely rational case $\r\bar\r$ is a true 
representation and we can express it explicitly. Here, the structure is 
better understood by dealing with the case of an arbitrary $n$-fold 
tensor product.

\begin{corollary}\label{canform}
Suppose $\A$ is completely rational, 
$\D=\A\otimes\cdots\otimes\A$ ($n$-fold 
tensor product) and let $\r$ be a soliton endomorphism equivalent 
to $\pi$ (see also Sect. \ref{tts}).

Then $[\r\bar\r] = [\bar\r\r]$ is a DHR sector, and we have the equality (as sectors)
\begin{equation}\label{multican}
\bar\r\r = \bigoplus_{i_0,i_1,\dots i_{n-1}} N^0_{i_0,i_1,\dots i_{n-1}}
\r_{i_0}\otimes\r_{i_1}\otimes\cdots\otimes\r_{i_{n-1}} \ ,
\end{equation}
where $N^0_{i_0,i_1,\dots i_{n-1}}$ is the multiplicity of the 
identity sector in the product $\r_{i_0}\cdot\r_{i_1}\cdots\r_{i_{n-1}}$
and the sum is taken over all irreducible sectors of $\A$. 
\end{corollary}
\proof
Formula (\ref{multican}) for $\bar\r\r$ follows immediately by Prop. 
\ref{intercan}, which gives $\bar\r\r$ in terms of the 
formula for the canonical endomorphism of the $n$-interval inclusion given
in \cite{KLM}  in the completely rational case.

To show that $\bar\r\r$ is equivalent to $\r\bar\r$ note that by 
eq. (\ref{conj}) we have, setting $\b\equiv\b_p=\b^2$,
\[
\r\bar\r= \r\b\r\b = \b(\b\r\b\r)\b =\b(\bar\r\r)\b\ ,
\]
that, combined with formula (\ref{multican}) gives
\begin{multline}
\r\bar\r 
=\bigoplus_{i_0,i_1,\dots i_{n-1}} N^0_{i_0,i_1,\dots i_{n-1}}
\r_{i_{p(0)}}\otimes\r_{i_{p(1)}}\otimes\cdots\otimes\r_{i_{p(n-1)}}\\
= \bigoplus_{i_0,i_1,\dots i_{n-1}}
N^0_{i_{p^{-1}(0)},i_{p^{-1}(1)},\dots i_{p^{-1}(n-1)}}
\r_{i_0}\otimes\r_{i_1}\otimes\cdots\otimes\r_{i_{n-1}}
\end{multline}
which coincides with formula (\ref{multican}) because the $\r_i$'s 
form a commuting family.
\endproof
Note in particular the special case $n=2$ in Cor. \ref{canform} gives 
the formula
\[
\r^2=\bar\r\r = \bigoplus_i \r_i\otimes\bar\r_i \ ,
\]
\section{Split \& $\mu_\A <\infty$ imply strong additivity}
\label{stradd}
Before deriving the strong additivity property from the finite $\mu$-index 
assumption, we recall some basic facts about simple subfactors 
\cite{L0}. Let $M$ 
be a factor in a standard form on a Hilbert space $\H$ with modular 
conjugation $J$. A subfactor $N\subset M$ is \emph{simple} if
\[
N\vee JNJ =B(\H)\ .
\]
In other words $N$ is a simple subfactor iff $N'\cap M_1=\mathbb C$ 
where $M_1\equiv JN'J$ is the basic extension in the sense of Jones 
\cite{J}; in particular $N'\cap M=\mathbb C$. 

If $N$ is a simple subfactor and there exists a normal conditional 
expectation $\e$ from $M$ onto $N$, then $N=M$. Indeed the 
expectation 
is faithful and the Takesaki-Jones projection implementing $\e$ 
belongs to $N'\cap M_1=\mathbb C$, thus $\e$ is the identity. In 
particular
\[
N\subset M\; \text{simple}\;\ \&\;\ [M:N]<\infty \Rightarrow N=M\ ,
\]
which is the implication we are going to use.
\smallskip

We now return to a local conformal net $\A$. 
We shall denote by $\A^d$ the dual net of $\A$ on $\mathbb R$, namely
$\A^d(I)\equiv\A(\mathbb R\setminus I)'$.
\begin{lemma}\label{ab1}
Let $\A$ be a local M\"obius covariant net. If $I\subset 
\mathbb R$ is a bounded interval and $I_1,I_2$ the intervals obtained 
by removing a point from $I$, we have:
\begin{itemize}
\item[$(a)$] $\A(I_1)\vee\A(I_2)\subset\A(I)\subset\A^d(I)$ is a 
basic extension. In particular 
\[
[\A^d(I):\A(I)]=[\A(I):\A(I_1)\vee\A(I_2)].
\]
\item[$(b)$] $\A(I_1)\vee\A(I_2)\subset\A(I)$ is a simple subfactor 
$\Leftrightarrow$ $\A(I_1)\vee\A(I_2)\subset\A^d(I)$ is irreducible 
$\Leftrightarrow$ $\A$ is 4-regular.
\end{itemize}
\end{lemma}
\proof
$(a)$: By dilation-translation covariance we can assume that 
$I_1=(-1,0)$, 
$I_2=(0,1)$, $I=(-1,1)$. The modular conjugation $J$ of $M\equiv 
\A(-1,1)$ is 
associated with the ray inversion map $t\to - 1/t$. With $N= 
\A(-1,0)\vee\A(0,1)$ we then have:
\begin{multline}
M_1\equiv JN'J=J\big (\A(-1,0)\vee\A(0,1)\big )'J
= \big (\A(-\infty,-1)\vee\A(1,\infty)\big )' =\\
\A(-1,\infty)\cap\A(-\infty,1) =\A^d(-1,1)
\end{multline}
$(b)$: This follows because
\begin{multline}
N\vee JNJ=\big (\A(-1,0)\vee\A(0,1)\big )\vee J\big 
(\A(-1,0)\vee\A(0,1)\big )J\\
= \big (\A(-1,0)\vee\A(0,1)\big )\vee \big 
(\A(-\infty,0)\vee\A(0,\infty)\big )
\end{multline}
which is equal to $B(\H)$ iff $\A$ is 4-regular.
\endproof
\begin{lemma}\label{4regstradd}
Let $\A$ be a local M\"obius covariant net. If $I\subset 
\mathbb R$ is a bounded interval and $I_1,I_2$ the intervals obtained 
by removing a point from $I$. Assume $[\A^d(I):\A(I)]<\infty$. We have:
\[
 \text{\rm $\A$ is 4-regular} \Rightarrow  \text{\rm $\A$ is strongly 
additive.}
\]
\end{lemma}
\proof
If $\A$ is 4-regular then $\A(I_1)\vee\A(I_2)\subset\A(I)$ is a 
simple subfactor by Lemma \ref{ab1}. On the other hand there exists a 
normal 
expectation $\A(I)\to\A(I_1)\vee\A(I_2)$ by the finite index 
assumption and Lemma \ref{ab1}.  But there is no normal expectation 
onto a simple subfactor, unless the inclusion is trivial. 
Thus $\A(I_1)\vee\A(I_2)=\A(I)$, i.e. $\A$ is strongly additive.
\endproof
\begin{theorem}\label{comprat}
Let $\A$ be a local conformal net with the split property. If 
$\mu_\A$ 
is finite, then $\A$ is strongly additive (thus completely rational).
\end{theorem}
\proof
If $\mu_\A <\infty$, then the 2-interval inclusion is irreducible by 
Lemma \ref{irr2intincl}, hence $\A$ is 4-regular. By the following 
Lemma \ref{ee} and Lemma \ref{4regstradd} we get the thesis.
\endproof
Let $\A$ be a split local conformal net.  If $\mu_{\A}<\infty$ we shall denote by 
$\e_E:\hat\A(E)\to\A(E)$ the conditional expectation associated 
with the 2-interval $E$ (unique by Lemma \ref{irr2intincl}).  
The following lemma is contained in \cite{KLM}. 
\begin{lemma}\label{ee}
Assume that the $\mu$-index of $\A$ is finite. Given a bounded 
interval $I\in\I$, there is a finite index expectation 
$\e_{I}:\A^d(I)\to\A(I)$.
\end{lemma}
\proof
Consider a decreasing sequence of 2-intervals $E_n\equiv I\cup I_n$ 
where $-1\in 
I_n$ and $\cap_n I_n=\{-1\}$. As shown in \cite{KLM}
\[
\A(E_n)\searrow \A(I),\quad \hat\A(E_n)\searrow \A^d (I)\ .
\]
As in  Prop. 2 of \cite{KLM}, any weak limit point $\e_I$ of $\e_{E_n}\res 
\hat\A^d(I)$ (as a map $\A^d(I)\to\A(E_1)$) is a finite index 
expectation from $\A^d(I)$ to $\A(I)$. 
\endproof
\section{Topological sectors and an index theorem}
\label{tts}

In this section we generalize the previous construction to the case 
of 
cyclic orbifold based on  a local conformal net $\A$ with the split 
property. 

Let $\zeta$ be a point of $S^1$ and 
$h:S^1\setminus\{\zeta\}\simeq\mathbb R\to S^1$ a smooth injective 
map 
which is smooth also at $\pm\infty$, namely the left and right limits
$\lim_{z\to\zeta^{\pm}}\frac{{\rm d}^n h}{{\rm d}z^n}$ exist for all 
$n$.

The range $h(S^1\setminus\{\zeta\})$ is either $S^1$ minus a 
point or a (proper) interval of $S^1$.

With $I\in\I$, $\zeta\notin I$, we set
\[
\Phi_{h,I}^{(\z)}\equiv \Ad U(k)\ ,
\]
where $k\in\Diff(S^1)$ and $k(z)=h(z)$ for all $z\in I$ and $U$ is 
the projective unitary representation of $\Diff(S^1)$ associated 
with $\A$.

Then $\Phi_{h,I}^{(\z)}$ does not depend on the choice of 
$k\in\Diff(S^1)$ and 
\[
\Phi_{h}^{(\z)}:I\mapsto \Phi_{h,I}^{(\z)}
\]
is a well defined soliton of $\A_0\equiv\A\restriction\mathbb R$.

Clearly $\Phi_h^{(\z)}(\A_0(\mathbb R))''=\A(h(S^1\setminus\{\z\}))''$, 
thus $\Phi_h^{(\z)}$ is irreducible if the range of $h$ is  
dense, otherwise it is a type III factor representation.  It is easy 
to see that, in the last case, $\Phi_h^{(\z)}$ does not depend on $h$ 
up to unitary equivalence.

Let now $f:S^1\to S^1$ be a smooth, locally injective map of degree 
deg$f=n\geq 1$. Choosing $\z\in S^1$, there are $n$ right inverses 
$h_i$, $i=0,1,\dots n-1$, for $f$; namely there are $n$ injective 
smooth maps $h_i:S^1\setminus\{\z\}\to S^1$ such that $f(h_i(z))=z$, 
$z\in S^1\setminus\{\z\}$. The $h_i$'s are smooth also at $\pm\infty$.

Note that the ranges $h_i(S^1\setminus\{\z\})$ are $n$ pairwise 
disjoint intervals of $S^1$, thus we may fix the labels of the $h_i$'s so that 
these intervals are counterclockwise ordered, namely
we have $h_0(-\z)<h_1(-\z)<\dots<h_{n-1}(-\z)<h_0(-\z)$. 

Of course any other possible choice for the $h_i$'s is 
associated with an element $p$ of the permutation group $\mathbb P_n$ 
on $\mathbb Z_n$, namely we can consider the sequence $h_{p(0)},h_{p(1)},\dots$.

For any interval $I$ of $\mathbb R$, we set
\begin{equation}\label{ts}
\pi_{f,I}^{(\z)}\equiv \chi_I\cdot 
(\Phi^{(\z)}_{h_{0},I}\otimes\Phi^{(\z)}_{h_1,I}\otimes\cdots\otimes\Phi^{(\z)}_{h_{n-1},I})\ ,
\end{equation}
where $\chi_I$ is 
the natural isomorphism from $\A(I_0)\otimes\cdots\otimes\A(I_{n-1})$ to 
$\A(I_0)\vee\cdots\vee\A(I_{n-1})$ given by the split property, 
with $I_k\equiv h_k(I)$. 
Clearly $\pi^{(\z)}_f$ is a soliton of 
$\D_0\equiv\A_0\otimes\A_0\otimes\cdots\otimes\A_0$ ($n$-fold 
tensor product).

If we order the right inverses $h_i$'s according to the permutation 
$p$ as above, we shall denote the corresponding soliton by 
$\pi_{f,p}$, thus $\pi_f\equiv\pi_{f,\text{id}}$. Clearly
\[
\pi_{f,p}=\pi_f\cdot\b_p
\]
where $\b$ is the natural action of $\mathbb P_n$ on $\D$.
\begin{proposition}
\label{pi}
Fix $\z=-1$ and denote $\pi^{(\z)}_f$ simply by $\pi_f$.

$(a)$: If $f_0$ has ${\rm deg}f_0 = {\rm deg}f$, then $\pi_{f_0}$ is 
unitary equivalent to $\pi_{f,p}$ for some $p\in\mathbb P_n$.

$(b)$: $\pi_{f,p}$ depends only on ${\rm deg}f$ and $p$ up to unitary 
equivalence.

$(c)$: ${\rm Index}(\pi_f) =\mu_{\A}^{n-1}$.

$(d)$: The conjugate of $\pi_f$ is given by 
\[
\bar\pi_f = \pi_{\bar f,p}
\]
where $\bar f(z)\equiv \overline{f(\bar z)}$ and $p$ is the inverse
automorphism $m\mapsto -m$ of $\mathbb Z_n$.
\end{proposition}
\proof
$(a)$: If $f_0:S^1\to S^1$ is a an injective smooth map and 
$\text{deg}f_0 = \text{deg}f$, there exists a $h\in\Diff(S^1)$ such 
that 
$f_0 = f\cdot h$. Then the $h^{-1}\cdot h_i$ are right inverses for 
$f_0$ 
and we have $\Phi_{h^{-1}\cdot h_i}=\Ad U(h)^*\cdot \Phi_{h_{p(i)}}$ 
for some $p\in \mathbb P_n$, so 
$U(h)$ implements a unitary equivalence between $\pi^{(\z)}_{f,p}$ and 
$\pi^{(\z)}_{f_0}$.

$(b)$: This is clear from the proof of $(a)$.

$(c)$: An obvious extension of Prop. \ref{intercan} shows that the 
index of $\pi_f$ is equal to the index of the $n$-interval inclusion, 
therefore by \cite{KLM} we have ${\rm Index}(\pi_f) =\mu_{\A}^{n-1}$.

$(d)$: if $\r$ is a soliton endomorphism of $\D$ localized $S^+$,
the formula in \cite[Th. 4.1]{GL1} gives $\bar\r=j\cdot\r\cdot j$, 
where $j=\text{Ad}J$ with $J$ the modular conjugation of $(\D(S^+),\Omega)$.
As we are interested in $\bar\pi_f$ up to unitary equivalence, we 
then have
\[
\bar\pi_f = j_0\cdot\pi_f\cdot j
\]
where $j_0\equiv\text{Ad}J_0$ with $J_0$
any unitary involution on the Hilbert space $\H$ of $\A$. Let $J_0$ 
then be the modular conjugation of $(\A(S^+),\Omega)$, thus 
$j=j_0\otimes\cdots\otimes j_0$. With the above notations let $x_0,x_1 , 
\dots x_{n-1}\in\A(\bar I)$ where $\bar I$ denotes here the conjugate 
interval of $I\subset S^1\setminus\{-1\}$. We have
\begin{align*}
\bar\pi_{f,\bar I}(x_0\otimes\cdots \otimes x_{n-1}) 
&= j_0(\pi_{f,\bar I}(j(x_0\otimes\cdots \otimes x_{n-1}))\\
&= j_0(\pi_{f,\bar I}(j_0(x_0)\otimes\cdots \otimes j_0(x_{n-1}))\\
&= j_0(\Phi^{(\z)}_{h_0,\bar I}(j_0(x_0))\cdots 
\Phi^{(\z)}_{h_{n-1},\bar I}(j_0( x_{n-1}))\\
&= \chi_{I}\cdot (j_0\cdot\Phi^{(\z)}_{h_0,\bar I}\cdot 
j_0\otimes\cdots\otimes 
j_0\cdot\Phi^{(\z)}_{h_{n-1},\bar I}\cdot j_0)(x_0\otimes\cdots 
\otimes x_{n-1})\\
&= \chi_{I}\cdot(\Phi^{(\z)}_{\bar h_{p(0)}, I}\otimes\cdots 
\otimes\Phi^{(\z)}_{\bar h_{p(n-1)}, I})(x_0\otimes\cdots \otimes x_{n-1})\\
&=\pi_{\bar f,I}(x_0\otimes\cdots \otimes x_{n-1}) .
\end{align*}
Here $p\in \mathbb P_n$ is the re-labeling of the right inverses $\bar h_i$ 
of $\bar f$ associated with the map $z\mapsto \bar z$ on the circle. 
It can be checked immediately in the case $f(z)=z^n$ that $p(k)=n-k$. 
\endproof
We shall now see the sector $[\pi_f]$ is independent of the choice of the
initial interval in the counterclockwise order associated with the 
$h_i$'s. Thus $[\pi]$ and $[\bar\pi_f]$ are the unique sectors associated 
respectively with any counterclockwise/clockwise ordering of the $h_i$'s. 
\begin{proposition}
\label{pi0}
$(a)$: If  $p\in\mathbb P_n$ is a cyclic permutation, then $\pi_f$ is
unitarily equivalent to $\pi_{f,p}$.

$(b)$: $\pi_f$ is irreducible if and only if $\A$ is $n$-regular.
\end{proposition}
\proof
$(a)$: It suffices to consider the case $f(z)=z^n$. With the choice 
of the $n$-th root function $\sqrt[n]{z}$ with  discontinuity 
at $-1$, we may order 
counterclockwise the right inverses by setting $h_{\ell}\equiv 
e^{\frac{2\pi \ell i}{n}}h_0$, $\ell\in 0,1,\dots n-1$.

Thus for any $j\in\mathbb Z_n$, 
$h_{\ell + j}=R^j\cdot h_{\ell}$, for all $\ell\in\mathbb 
Z_n$, where $R\equiv R(\frac{2\pi}{n})$ denotes the 
rotation on $S^1$ of angle $\frac{2\pi}{n}$, and so $U(h_{\ell + j})=
U(R^j)U(h_{\ell})$ (up to a phase factor).

If $p$ is the cyclic permutation $\ell\mapsto \ell + j$ on $\mathbb 
Z_n$, it follows that
\begin{align}
\pi_{f,p,I} & = \chi_I\cdot
(\Phi^{(\z)}_{h_{j},I}\otimes\Phi^{(\z)}_{h_{j+1},I}
\otimes\cdots\otimes\Phi^{(\z)}_{h_{j+n-1},I}) \\
& = \chi_I\cdot 
(\Ad U(R^j)\cdot\Phi^{(\z)}_{h_{0},I}
\otimes\Ad U(R^j)\cdot\Phi^{(\z)}_{h_1,I}
\otimes\cdots\otimes\Ad U(R^j)\cdot\Phi^{(\z)}_{h_{n-1},I})\\
& = \Ad U(R^j)\cdot\chi_I\cdot
(\Phi^{(\z)}_{h_{0},I}\otimes\Phi^{(\z)}_{h_{1},I}
\otimes\cdots\otimes\Phi^{(\z)}_{h_{n-1},I}) \\
& = \Ad U(R^j)\cdot \pi_{f,I}
\end{align}

$(b)$: As $I$ varies in the intervals of $S^1\setminus\{-1\}$, 
$\pi_I(\D(I))=\A(I_0)\vee\cdots\A(I_{n-1})$ generates 
$\A(S^1\setminus F)$ where $F$ is the set of $n$ points obtained by 
removing $\cup_i h_{i}(S^1\setminus\{-1\})$ from $S^1$, hence the thesis.
\endproof
\noindent
{\bf Remark.}
As already said, $2$-regularity is automatic for any M\"{o}bius 
covariant local net; but there are examples of M\"{o}bius 
covariant local nets that are not 3-regular \cite{GLW}. We conjecture 
that every diffeomorphism covariant local net is 
automatically $n$-regular for any $n$.
\medskip

As $\z$ varies, the $\Phi^{(\z)}_{k}$'s undergo permutations among 
them, indeed cyclic permutations that, with a proper labeling, 
correspond to the cyclic permutations on $(0,1,\dots, n-1)$. The 
restriction
\[
\t_f\equiv \pi^{(\z)}_f\restriction 
(\A\otimes\A\cdots\otimes\A)^{\mathbb Z_n}
\]
is therefore a DHR representation of
$(\A\otimes\A\cdots\otimes\A)^{\mathbb Z_n}$, independent of $\z$ up 
to unitary equivalence.

In the following we shall denote by 
$\I^{(n)}$ the set of all $n$-intervals of $S^1$, not necessarily 
symmetric (union of $n$ intervals with pairwise disjoint closures).

\begin{theorem}
\label{indexformula}

$(a)$: $\t_f$ depends only on $n={\rm deg}f$ up to unitary 
equivalence.

$(b)$: $\t_f$ is diffeomorphism covariant; the corresponding 
projective unitary 
representation of $\Diff^{(\infty)}(S^1)$ is unitary equivalent 
to the projective unitary representation $U^{(n)}=U\cdot M^{(n)}$ 
of $\Diff^{(n)}(S^1)$.

$(c)$: The following formula for the index holds:
\[
{\rm Index}(\t_f) =n^2 \mu_{\A}^{n-1}\ .
\]

$(d)$: $\t_f$ is direct sum of $n$ diffeomorphism 
covariant representations 
$\t^{(0)}_{f}, \t^{(1)}_{f},\dots,\t^{(n-1)}_{f} $ of 
$(\A\otimes\A\cdots\otimes\A)^{\mathbb Z_n}$. Each 
$\t^{(i)}_{f} $ is irreducible.

$(e)$ We may choose our labels  
so that, for every $i=0,1,\dots, n-1$,
\begin{gather*}
{\rm spin}(\t^{(i)}_{f}) = \frac{i}{n} + \frac{n^2 - 1}{24 n}c\ , \\
{\rm Index}(\t^{(i)}_{f}) = \mu_{\A}^{n-1}\ ,  
\end{gather*}
where, in the last equation, we assume $\mu_{\A}<\infty$.
\end{theorem}
\proof
$(a)$: Immediate by $(a)$ of Prop. \ref{pi}.

$(b)$: Because of the above point, it suffices to consider the case 
$f(z)=z^n$. Then the covariance follows by the characterization of 
the 
map $M^{(n)}$ in Prop. \ref{cover2} expressed by the commutativity of 
the diagram (\ref{cdn}).

$(c)$: Analogously as in Proposition \ref{iso2}, the inclusion 
$\pi_f(M)\subset\pi_f(M')'$, 
$M=(\A\otimes\cdots\otimes\A)(0,\infty)$, is 
isomorphic to the $n$-interval inclusion $\A(E)\subset\hat\A(E)$, 
$E\in\I^{(n)}$. If $\mu_{\A}<\infty$, then $\A$ is completely 
rational and the index formula in \cite{KLM} gives
\[
\text{Index}(\pi_f)= [\hat\A(E):\A(E)] =\mu_{\A}^{n-1} \ .
\]
As $\t_f$ is the restriction of $\pi^{(\z)}_f$ to a $n$-index subnet 
we then have
\[
\text{Index}(\t_f)
=[(\A\otimes\cdots\otimes\A)^{\mathbb Z_n}:\A\otimes\cdots\otimes\A]\cdot
\text{Index}(\pi_f)=n^2\mu_{\A}^{n-1}\ .
\]

$(d)$: 
Fix an interval $I_0$ and a 
unitary $v\in\D(I_0)$, $v^n =1$, that implements the  action on $\B$ 
dual to cyclic permutations. Then $\D(I)=\{\B(I),v\}''$ for all 
intervals $I\supset I_0$, hence
\[
\{\bigvee_{\z \notin\bar I}\t_{f,I}(\B(I)),\pi_f^{(\z)}(v)\}''=
\big(\bigvee_{\z \notin\bar I}\pi_{f,I}^{(\z)}(\D(I))\big )''=\A(S^1\setminus F)\ ,
\]
where $F$ is an $n$-point subset of $S^1$ (the complement of 
$\cup_i h_i(S^1\setminus \{\z\})$), that depends on $\z$. 

Now $\t_f$ is a DHR representation, so we may vary the point $\z$ and 
get
\[
\{\bigvee_{ I\in\I}\t_{f,I}(\B(I)),\pi_f^{(\z)}(v)\}''=
\bigvee_{I\in\I}\pi_{f,I}^{(\z)}(\D(I))=\A(S^1)=B(\H)\ ,
\]
where $\z\notin\bar I$ varies with $I$. As $\pi_f^{(\z)}(v)$ normalizes
$\bigvee_{ I\in\I}\t_{f,I}(\B(I))$, it follows as in Cor. \ref{tirr} 
for the case $n=2$ that the latter is the commutant of $\pi_f(v)$ and 
$\t_f$ has exactly $n$ irreducible components.

$(e)$: As in the case $n=2$, the covariance of $\t_f$ is given by a 
unitary representation of $\Diff(S^1)$ equivalent to $U^{(n)}=U\cdot 
M^{(n)}$. Thus the conformal Hamiltonian $L'_0$ in the 
representation $\t_f$ is unitarily equivalent to the one given by 
formula (\ref{spin}), and this readily implies that the spin of the 
$\t^{(i)}_{f}$'s are as stated, by a suitable choice of the index 
labels. We will have additional information on these labels in Sect. 8.3 after
(\ref{disa}). 
Concerning the formula for the index, by (\ref{at}) we have
$d(\tau^{(i)}_{f})=d(\pi_f).$
By point  $(c)$ we have $ d(\pi_f) = \sqrt{\mu_\A^{n-1}}$, thus 
${\rm Index}(\tau^{(i)}_{f})=\mu_\A^{n-1}$. 
\endproof
\subsection{Extension to non-vacuum representation case}
The construction given above in Sect. \ref{tts} extends to the case 
where one replaces the vacuum representation with another covariant 
representation $\lambda$ (cf. \cite{GL1,Ks} and Appendix \ref{covar} for the 
covariance condition). This extension generates new sectors and 
will be later used. Here we merely outline the construction, but all the above 
results have natural extensions in this setting.

Let $\lambda$ be a covariant representation of $\A$. Given an 
interval $I\subset S^1\setminus\{\z\}$,  we set
\[
\pi_{f_\lambda,I}(x)= \lambda_J(\pi_{f,I}(x))\ 
,\quad x\in\D(I)\ ,
\]
where $\pi_{f,I}\equiv\pi_{f,I}^{(\z)}$ is defined as in (\ref{ts}), 
and $J$ is any interval which contains $I_0\cup I_1\cup...\cup I_{n-1}$. 
\begin{proposition}\label{non-v}
The above definition is independent of the choice of $J$, thus 
$\pi_{f_\lambda,I}$ is a well defined soliton of $\D$. 

We can choose an interval $I$ with $\z$ as a boundary point
of $I$ such that  $\pi_f$ , $\pi_{f_\lambda}$ and $\lambda$ are localized on
$I$. Denote by $\tilde \pi_f, \tilde \pi_{f_\lambda}$ and 
$(\l,1,1,...,1):=\l\otimes \iota \otimes \iota\cdots \otimes 
\iota\restriction\D(I)$ respectively the corresponding endomorphisms of  $\D(I)$. 
Then  as sectors of $\D(I)$ we have
\[
[\tilde\pi_{f_\lambda}]= 
[\tilde\pi_f\cdot(\l,1,1,...1)\ ].
\]
\end{proposition}
\proof
If $J_1$ is another interval which contains $I_0\cup I_1\cup\cdots\cup I_{n-1}$, 
we need to show that $ \pi_{\lambda,J_1}(x)=  \pi_{\lambda,J}(x),
\forall x\in  \A(I)\otimes\cdots\otimes\A(I)$. It is sufficient to prove
this for $x=x_0\otimes\cdots\otimes x_{n-1}$, $x_i\in  \A(I)$, $0\leq i\leq n-1$. 
By isotony, we have 
\begin{multline*}
\pi_{f_\lambda,J_1}(x_0\otimes\cdots\otimes x_{n-1}) 
= \l_{J_1}(\Phi_{h_0,I}(x_0))\cdots
\l_{J_1}(\Phi_{h_{n-1},I}(x_{n-1}))\\
= \l_{J}(\Phi_{h_0,I}(x_0))\cdots
\l_{J}(\Phi_{h_{n-1},I}(x_{n-1}))
= \pi_{f_\lambda,J}(x_0\otimes\cdots\otimes x_{n-1} )\ . 
\end{multline*}
This shows that the above definition is independent of the choice of $J$.
\par
As for the last formula, we may assume that $\z=-1$, 
$I=S^+$ (the upper half circle), 
$f(z)=z^n$, $h_0$ 
is the $n^{\text{th}}$-root function on 
$I$ with $h_0(1)=1$ so that $I_0 \subset I$, and $h_0 \in  
\Diff(J_0)$ for some interval 
$J_0 \supset I$, i.e. $h_0\in\Diff(S^1)$ and $h_0$ acts identically on $J'_0$.  
We may further assume that $\l$ is localized in $I_0$.
By our assumption $U(h_0)\in \A(J_0)$, and we claim
that
\begin{equation}\label{covariance}
\lambda_{J_0} (U(h_0) \lambda_{J_1}(x) \lambda_{J_0} (U(h_0))^*
=  \lambda_{h_0(J_1)}(U(h_0)x  U(h_0)^* ),\ 
\forall x\in \A(J_1), \forall J_1\in\I\ .
\end{equation}
This can be checked as follows: If $\bar J_0\cup \bar J_1\neq S^1$, 
then we can 
find an interval $J_2$ such that $J_0\cup J_1\subset  J_2$, and in this
case
\[
\lambda_{J_0} (U(h_0) \lambda_{J_1}(x) \lambda_{J_0} (U(h_0))^*
=  \lambda_{J_2}(U(h_0)x  U(h_0)^* ), \ \forall x\in \A(J_1)\ ;
\]
note that $U(h_0)x  U(h_0)^*\in \A(h_0(J_1))$, and
$h_0(J_1) \subset J_2$ , so by isotony we have 
\[
\lambda_{J_2}(U(h_0)x  U(h_0)^* )=  \lambda_{h_0(J_1)}(U(h_0)x  
U(h_0)^* )\ .
\]
In general we cover $J_1$ by a set of sub-intervals $J_k\subset J_1, 
2\leq k\leq m$ 
such that
$\bar J_k \cup \bar J_0 \neq S^1$. By additivity of conformal nets we 
have $\A (J_1) =\vee_{2\leq k\leq m} \A(J_k)$, and since the 
equation (\ref{covariance}) is true for any $x\in  \A(J_k), 2\leq k\leq m$,
it follows that we have proved  equation (\ref{covariance}). \par
Define $z_\lambda(h_0):= \lambda_{J_0} (U(h_0)) U(h_0)^*$. From
(\ref{covariance}) we have
$$
\lambda_{h_0(J_1)}(\Ad U({h_0})(x)) = 
z_\lambda(h_0) \Ad U({h_0})(\lambda_{J_1}(x)) z_\lambda(h_0)^*
, \forall x\in \A(J_1), \forall J_1.$$
Set $J_1=I'$, we conclude from the above equation that 
$z_\lambda(h_0)\in \A(I_0')'= \A(I_0)$. It follows that
for all $x_0\otimes x_1\cdots\otimes x_{n-1} \in \D(I)$, we have
\begin{multline*}
\lambda_{I_0} (\Ad U({h_0})(x_0))\otimes 
\Ad U({h_1})(x_1)\otimes\cdots\otimes \Ad U({h_{n-1}})(x_{n-1})
\\
= \Ad{z_\lambda (h_0)}\cdot\Ad U({h_0}(\l_I(x_0)) \otimes  
\Ad U({h_1})(x_1)\otimes\cdots\otimes \Ad U({h_{n-1}})(x_{n-1}))
\end{multline*}
where $h_1,...,h_{n-1}$ are defined as in (\ref{ts}). 
Therefore on $ \D(I)$
\begin{equation}\label{cov2}
\pi_{f_\lambda,I} =  \Ad{z_\lambda(h_0)}\cdot \pi_{f,I}\cdot (\l,1,1,...1).
\end{equation}
Let $U_{I'}: \H\rightarrow \H\otimes \H\otimes\cdots
\otimes \H$ ($n$-tensor factors)
be a unitary such that
$ U_{I'} \pi_{f,I'}(\cdot) U_{I'}^* = id $ on $\D(I')$. Then both 
$ \tilde  \pi_f:= U_{I'}\pi_{f,I}(\cdot) U_{I'}^*$ and 
$\tilde  \pi_{f_\lambda}:=U_{I'}\pi_{f_\lambda,I}(\cdot)U_{I'}^*$
are endomorphisms of $\D(I)$, and we have
$   \tilde  \pi_{f_\lambda} =  \Ad {U_{I'} z_\lambda(h_0)  U_{I'}^*}
\cdot\tilde  \pi_f \cdot (\l, 1,1,...,1)$ by (\ref{cov2}).  
Therefore, as sectors of $\D(I)$, we have
$$
[  \tilde  \pi_{f_\lambda}  ]=[ \tilde  \pi_f \cdot (
\l, 1,1,...,1)] 
$$
since $U_{I'} z_\lambda(h_0)  U_{I'}^*\in U_{I'} \A(I_0) U_{I'}^*
\subset \D(I')'=\D(I)$.
\endproof
\section{Some consequences}
We now discuss a few consequences of our results. The first two ones 
follow immediately from the implication ``rationality $\Rightarrow$ 
complete rationality'' because of the corresponding results in 
\cite{KLM} in the 
completely rational case.
\subsection{Rationality implies modularity}
The first consequence concerns the invertibility of the matrices $T$ 
and $S$ in a rational model, see \cite{R2}. 
This property has long been expected and 
is at the basis most analysis, in particular concerning Topological 
QFT, cf. for example \cite{GW}. 

We shall say that a local conformal net $\A$ is \emph{rational} if 
there are only finitely many irreducible sectors and all of them 
have a conjugate sector, i.e. they have finite index \cite{L1,GL2}. 
Assuming the split property, then every sector is direct sum 
of irreducible sectors, cf. \cite{KLM}.

In the paper \cite{KLM} the modularity has been proved for a 
completely rational local M\"{o}bius covariant net. By our results, complete 
rationality is equivalent to rationality for a local conformal net 
with the split property. Hence we have: 
\begin{theorem}
Let $\A$ be a conformal net with the split property.  
If $\A$ is rational, then the tensor category of representations of 
$\A$ is modular, i.e. the braiding symmetry is non-degenerate.
\end{theorem}
\subsection{The $\mu$-index is always equal to the global index}
The equality of the $\mu$-index with the global index has been proved 
in \cite{KLM} in the completely rational case. The extension of
this equality to the case of infinite $\mu$-index is not covered by 
that work, in particular there was no argument to show that if there 
is no non-trivial sector then Haag duality holds for multi-connected 
regions. This is given here below. 
\begin{theorem}
Let $\A$ be a conformal net with the split property.  Then
\begin{equation}\label{globindex}
\mu_\A =\sum_i d(\r_i)^2\ ,
\end{equation}
where the sum is taken over all irreducible sectors or, equivalently, 
over all the irreducible sectors that are diffeomorphism covariant 
with positive energy.
\end{theorem}
\proof
If $\mu_{\A} <\infty$ then $\A$ is completely rational by Theorem 
\ref{comprat}, thus the formula holds by \cite{KLM}.

If $\mu_{\A}=\infty$ either there exists an irreducible sector with 
infinite index and formula (\ref{globindex}) obviously holds, or by 
Th. \ref{comprat} there are (uncountably) infinitely many irreducible 
sectors, thus (\ref{globindex}) holds because $d(\r_i)\geq 1$.
\endproof
\begin{corollary}
Let $\A$ be a conformal net with the split property. The following 
are equivalent:
\begin{itemize}
\item [$(i)$] $\A$ has no non-trivial representation,
\item [$(ii)$] Haag duality holds for some $n$-intervals $E$ for 
some $n\geq 2$: $\A(E)'= \A(E')$,
\item [$(iii)$] Haag duality holds for all $n$-intervals: $\A(E)'= 
\A(E')$ for all $E\in\I^{(n)}$, $\forall n\in\mathbb N$.
\end{itemize} 
\end{corollary}
\proof
By eq. (\ref{globindex}), $(i)$ holds iff $\mu_{\A}=1$, namely iff
$(ii)$ holds with $n=2$. In this case $\A$ is completely rational by Th. 
\ref{globindex} and the formula $[\hat\A(E):\A(E)]=\mu^{n-1}_{\A}$, 
$E\in\I^{(n)}$, in \cite{KLM} shows that also $(iii)$ holds.

It remains to show that $(ii)\Rightarrow (iii)$. Assume that $\A(E)'= 
\A(E')$ for some $n$-interval $E$, $n\geq 2$. Then $\A(E)'= 
\A(E')$ for all $n$-intervals $E$ by diffeomorphism covariance. Fix 
$E\in\I^{(n)}$ and $I$ one of its connected components. By considering
a decreasing sequence of intervals $I\supset 
I_1\supset I_2\supset\cdots$ shrinking to a point, it is rather 
immediate to check, by the split property, that Haag duality $\A(E)'= 
\A(E')$ holds for $n-1$-intervals. By iteration we get Haag duality 
for a $2$-interval and then conclude our proof as above.
\endproof
\subsection{Sectors with infinite statistics}
General properties of sectors with infinite dimension were studied in 
\cite{GL1} (see also \cite{BCL}), yet first examples have been constructed 
by Fredenhagen in \cite{F}, see below. 
A natural family of infinite dimensional irreducible sectors  
has recently been pointed out by Carpi \cite{C} in the Virasoro nets with 
$c>1$, following a conjecture by Rehren \cite{R1}.

The following theorem gives a natural and general
construction of irreducible sectors with infinite dimension, as a 
consequence of the index formula in Sect. \ref{tts}.
\begin{theorem}\label{infstat}
Let $\A$ be a conformal net with the split property. The following 
are equivalent:
\begin{itemize}
\item[$(i)$] $\A$ is not completely rational;
\item[$(ii)$] $(\A\otimes\A)^{\text{\rm flip}}$ has an irreducible 
sector with infinite dimension;
\item[$(iii)$] $(\A\otimes\cdots\otimes\A)^{\mathbb Z_n}$ has an 
irreducible sector with infinite dimension and diffeomorphism 
covariant with positive energy, any $n\geq 2$.
\end{itemize}
\end{theorem}
\proof
Clearly $(ii)$ or $(iii)$ imply that $\A$ is not completely rational 
(complete rationality if hereditary for finite-index subnets 
\cite{Xu4,L4}). On the other hand, if $\A$ is not completely 
rational, the topological sector $\tau_f$ of the cyclic $n$-orbifold
has infinite index by the index formula in Th. \ref{indexformula}. So 
one of the $n$ direct summands $\t^{(i)}_f$  
must have infinite index.  
\endproof
\subsubsection{Example}
Let $\A$ be the local conformal net on $S^1$ associated with the 
$U(1)$-current algebra. In the real line picture $\A$ is given by
\[
\A(I)\equiv \{W(f): f\in C_{\mathbb R}^{\infty}(\mathbb R),\ 
\text{supp}f \subset I\}''
\]
where $W$ is the representation of the Weyl commutation relations
\[
W(f)W(g)=e^{-i\int fg'}W(f+g)
\]
associated with the vacuum state $\o$
\[
\o(W(f))\equiv e^{-||f||^2},\quad ||f||^2\equiv\int_0^{\infty}|\tilde 
f(p)|^2p\text{d}p
\]
where $\tilde f$ is the Fourier transform of $f$.

The superselection structure of $\A$ is completely described in 
\cite{BMT}. There is a one parameter family $\{\a_q, q\in \mathbb 
R\}$ of irreducible sectors 
and all have index $1$. We can choose a representative of $\a_q$ as
\[
\a_q(W(f))\equiv e^{2i\int\! Ff}W(f), \quad F\in C^{\infty},\quad 
\int 
F = q\ .
\]
Now consider $\A\otimes\A$. By the argument in \cite{KLM} all 
irreducible sectors of $\A\otimes\A$ are tensor product sectors, 
namely have the form $\a_q\otimes\a_{q'}$, in particular they have 
index 
$1$.

Yet, the index $2$ subnet $(\A\otimes\A)^{\text{flip}}$ has an 
irreducible sector with infinite index, by Th. \ref{infstat} because 
$\A$ is not 
completely rational.

Fredenhagen  \cite{F} had shown that the subnet 
$(\A\otimes\A)^{SO(2)}\subset\A\otimes\A$ admits an 
infinite dimensional irreducible sector. 
In his case the subnet $(\A\otimes\A)^{SO(2)}\subset \A\otimes\A$ 
has infinite index. 
\section{Topological twisted sectors in the completely rational case}
\label{fus}

In this section we assume that $\D$ is a completely rational conformal 
net and $\B$ is the fixed point subnet of $\B$ under the proper action
of ${\mathbb Z_n}$ on $\D$ (cf. 2 of \cite{Xu3}).  

We note that we will be interested in the special case when
$\D:= \A\otimes\A...\otimes\A$ ($n$-fold tensor product) and  
$\B:=(\A\otimes \A...\otimes\A)^{\mathbb Z_n}$ the fixed point subnet 
of $\D$ under the action of cyclic permutations in Sect. \ref{fus}.

By Th. 2.9 of \cite{Xu3}, $\B$ is completely rational with
$\mu_{\B}= n^2 \mu_{\D}$. So  $\B$ has finitely many inequivalent irreducible
representations and the question is how to construct these 
representations
from those of $\D$. This question 
can be raised for the case of a general orbifold. 
An answer to this question is given in an example
of $\mathbb Z_2$ orbifold of a lattice by identifying the orbifold
with a coset whose irreducible representations are known (cf. Sect. 3 of
\cite{Xu3}). Partially motivated by this question for the case of
cyclic permutations, \cite{Xu6} and 
\cite{BDM}, 
we were led to the constructions of Sect. \ref{first} and \ref{tts}. We will see that  
the topological construction of  Sect. \ref{tts} and its generalizations lead to 
a satisfying answer to the question for $n=2,3,4$ and plays an 
important 
role in the general description of cyclic orbifold.\par
In this section we will make use of computations of sectors 
extensively as
in \cite{Xu1}. Let us first recall some preliminaries about sectors. 
See
\cite{L1}, \cite{L1'} and \cite{L5} for more details.  
Let $M$ be an infinite factor
and  $\text{\rm End}(M)$ the semigroup of
 unit preserving endomorphisms of $M$.  
Let $\text{\rm Sect}(M)$ denote the quotient of $\text{\rm End}(M)$ 
modulo
unitary equivalence in $M$. We  denote by $[\rho]$ the image of
$\rho \in \text{\rm End}(M)$ in  $\text{\rm Sect}(M)$.\par
 It follows from
\cite{L1'} that $\text{\rm Sect}(M)$ is endowed
with a natural involution $\theta \rightarrow \bar \theta $  ;
moreover,  $\text{\rm Sect}(M)$ is a semiring.\par 

Let $\rho \in \text{\rm End}(M)$ and $\e$ be a normal
faithful conditional expectation
$\e:
M\rightarrow \rho(M)$.  We define a number $d_\e\geq 1$ (possibly
$\infty$) by:
$$
d_\e^{-2} :=\text{\rm Max} \{t\in [0, +\infty)|
\e (m_+) \geq t m_+, \forall m_+ \in M_+
\}$$ (Pimsner-Popa inequality in \cite{PP}).\par
 We define
$$
d(\rho) = \text{\rm Min}_\e \{ d_\e \},
$$ 
where the minimum is taken over $\e$ with $d_\e < \infty$ (otherwise 
we put $d(\r) =\infty$).  $d(\rho)$ is called the  dimension of  $\rho$. It is 
clear from the definition that  the dimension  of  $\rho$ 
depends only the sector $[\rho]$.

The properties of the dimension can be found in
\cite{L1}, \cite{L1'} and  \cite{L5}, see also \cite{Ko}.
We recall that $d(\r)<\infty$ is equivalent to the existence of a 
conjugate sector. \par
For $\lambda $, $\mu \in \text{\rm End}(M)$, let
$\text{\rm Hom}(\lambda , \mu )$ denote the space of intertwiners from
$\lambda $ to $\mu $, i.e. $a\in \text{\rm Hom}(\lambda , \mu )$ iff
$a \lambda (x) = \mu (x) a $ for any $x \in M$.
Assuming the dimension of $\l$ and $\mu$ to be finite,
$\text{\rm Hom}(\lambda , \mu )$  is a finite dimensional vector
space and we use $\langle  \lambda , \mu \rangle$ to denote
the dimension of this space.  $\langle  \lambda , \mu \rangle$
depends
only on $[\lambda ]$ and $[\mu ]$. Moreover we have
$\langle \nu \lambda , \mu \rangle =
\langle \lambda , \bar \nu \mu \rangle $,
$\langle \nu \lambda , \mu \rangle
= \langle \nu , \mu \bar \lambda \rangle $ which follows from 
Frobenius
duality (see \cite{L2} ).  We will also use the following
notation: if $\mu $ is a subsector of $\lambda $, we will write as
$\mu \prec \lambda $  or $\lambda \succ \mu $.  A sector
is said to be irreducible if it has only one subsector. 
Usually we will use Greek letters to denote sectors, but we 
will denote the identity sector by $1$ if no confusion arises. \par
Fix an interval $J_0$. Let $\gamma: \D(J_0)\rightarrow \B(J_0)$ be the 
canonical endomorphism from $ \D(J_0)$ to $\B(J_0)$ and let $\gamma_{\B}:=
\gamma\res \B(J_0)$. Note $[\gamma]=[1]+[g]+...+[g^{n-1}]$ as sectors of
$\D(J_0)$ and $[\gamma_{\B}]=[1]+[\sigma]+...+[\sigma^{n-1}]$ as sectors of
$\B(J_0)$. Here $[g^i]$ denotes the sector of $\D(J_0)$ which is the 
automorphism 
induced by $g^i$ and $\sigma$ is a DHR representation of $\B$ with
$[\sigma^n]= [1]$ where $[1]$ denotes the identity sector. 
We note that the notation  $[g^i]$ is an exception to our rule
of using Greek letters to denote sectors. All the sectors considered
in the rest of  Sect. \ref{fus} will be sectors of  $ \D(J_0)$ or  $ \B(J_0)$ as
should be clear from their definitions. 
All DHR representations will be assumed to be localized
on $J_0$ and have finite statistical dimensions. For simplicity of notations,
for a DHR representation
$\sigma_0$ of $\D$ or $\B$ localized on $J_0$,
we will use the same notation $\sigma_0$ to denote its restriction to
$ \D(J_0)$ or  $ \B(J_0)$
and we will make no distinction
between local and global intertwiners (cf. Appendix A) 
for DHR representations localized on $J_0$ since they are the same by the
strong additivity of $\D$ and $\B$.   
\subsection{Non-twisted sectors in general case}
We will denote by 
$\lambda$ the irreducible DHR
representations of $\D$ and by $\lambda_{\B}$ 
its restriction to $\B$.  $\lambda_{\B}$ and its 
irreducible summands will be called {\it non-twisted} representations 
(relative to $\D$).  An irreducible DHR representation of $\B$ is {\it 
twisted} if it is not non-twisted.  Our goal in this section is to 
characterize the nature of non-twisted representations.\par
Let $\sigma_1$ be a DHR representation of $\B$ localized on $J_0$. Recall
from \S \ref{3.1} the definition of  $\alpha_{\sigma_1}$. When 
restricted to
$ \D(J_0)$,  $\alpha_{\sigma_1}$ is an endomorphism of $ \D(J_0)$
(cf. (1) of Th. 3.1 in \cite{Xu1} or Cor. 3.2 of \cite{BE}), 
and we use the same notation  $\alpha_{\sigma_1}$
to denote this endomorphism. For the rest of Sect. \ref{fus}, 
$\alpha_{\sigma_1}$ will always be understood as  the 
endomorphism of $\D(J_0)$. 
The following lemma which follows essentially from 
\cite{Xu1} (also cf. \cite{BE}) will be used repeatedly:
\begin{lemma}\label{key}
Let $\sigma_1, \sigma_2$ (resp. $\lambda,\mu$) 
be DHR representations of $\B$ (resp. $\D$) localized on $J_0$. Then:\par
(1) $[\bar\alpha_{\sigma_1}]=[\alpha_{\bar\sigma_1}]$ as sectors of
$\D(J_0)$ and $d({\alpha_{\sigma_1}})= d({\sigma_1})$; \par
(2) 
$
\langle \alpha_{\sigma_1}, \alpha_{\sigma_2} \rangle
= \langle \sigma_1\gamma_{\B}, \sigma_2 \rangle,$
$
\langle \alpha_{\sigma_1},  \lambda\rangle
= \langle \sigma_1,  \lambda_{\B}\rangle;$\par
(3) $
[g\alpha_{\sigma_1}]= [\alpha_{\sigma_1}g],$
$[\lambda\alpha_{\sigma_1}]= 
[\alpha_{\sigma_1} \lambda];$
\par 
(4)
$
\langle g^i\alpha_{\sigma_1},  g^j\alpha_{\sigma_2} \rangle
=\delta_{ij} \langle \alpha_{\sigma_1},  \alpha_{\sigma_2} \rangle$,
$\langle g^i\alpha_{\sigma_1},  g^j  \lambda \rangle
=\delta_{ij} \langle \alpha_{\sigma_1},   
\lambda\rangle$, 
$\langle g^i\mu,  g^j  \lambda \rangle
=\delta_{ij} \langle \mu,   
\lambda\rangle$, $0\leq i,j\leq n-1. $
\end{lemma}
\proof
$(1)$ follows from Cor. 3.5 of \cite{Xu1}, $(2)$ follows from Th. 3.3 of
\cite{Xu1}, $(3)$ follows from Th. 3.6 of \cite{Xu1}, and $(4)$ follows 
from
Lemma 3.5 of  \cite{Xu1}.
\endproof
Note that $\mathbb Z_n$ acts on  $\lambda$ 
naturally by $g\lambda g^{-1}$: this is a DHR representation of 
$\D$ localized on the fixed interval $J_0$ and whose restriction to
$\D(J_0)$ is simply $g\cdot \lambda \cdot g^{-1}$.
Assume that the stabilizer of such an action on $\lambda$ is 
generated by
$g^{n_1}$ with $n_1k_1=n$. Then:
\begin{proposition}\label{nong}
$\lambda_{\B}$ decomposes into $k_1$ different 
irreducible pieces denoted by 
$(\lambda; \sigma^i)$, $0\leq i\leq {k_1-1}$.
Moreover $ [\alpha_{(\lambda; \sigma^i)}]=
\bigoplus_{0\leq k\leq {n_1-1}}[g^k\lambda g^{-k} ]$, 
$d({(\lambda; \sigma^i)})= n_1 d(\lambda)$, 
and
if  $[(\lambda; \sigma^i)]= [(\mu; \sigma^j)]$
then  there exists an integer $l$ such that $\mu= 
g^l \lambda g^{-l}.$ 
\end{proposition}
\proof
Let $\rho_1$ be an endomorphism of $\D(J_0)$ such that 
$\rho_1(\D(J_0))= \B(J_0)$ and $\rho_1\bar\rho_1= \gamma$. 
By \cite{LR} as sectors of $\B(J_0)$ 
we have $[\lambda_{\B}]= 
[\gamma \lambda\res \B(J_0)]$, 
it follows that 
\[
\Hom(\lambda_{\B}, 
\lambda_{\B})_{\B(J_0)}
\simeq \Hom(\bar\rho_1 \lambda\rho_1, 
\bar\rho_1 \lambda\rho_1)_{\D(J_0)}
\]
By  Frobenius duality we have
\[
\langle \lambda_{\B}, 
\lambda_{\B}\rangle
=\langle \lambda, 
 \gamma\lambda\gamma\rangle  
\]
For $0\leq i,j\leq n-1$, note that $ g^j\lambda g^{-j}$ 
is a DHR representation of $\D$, and  
by (4) of lemma \ref{key} we have
$\langle \lambda, g^i \lambda g^{-j} \rangle
= \langle \lambda, g^{i-j} g^j\lambda g^{-j} \rangle
= \delta_{ij} \langle \lambda, g^j\lambda g^{-j} \rangle.$
It follows 
that $
\langle \lambda, 
 \gamma\lambda\gamma\rangle = k_1.$  

Notice that $[g\rho_1]=[\rho_1], [\bar\rho_1 g ]=[\bar\rho_1 g]$. 
If we set $ \nu_1=\bar\rho_1, \nu=g^{n_1}, 
\nu_2=\lambda\rho_1$, 
we have $[\nu_1\nu]=[\tau_1], [\nu \nu_2]=[\nu_2]$, 
and $\nu$ has order $k_1$. 
Now apply Lemma 2.1
of \cite{Xu5} 
where  $a,\nu,b$ of \cite{Xu5}  
correspond to our $\nu_1,\nu, \nu_2$ respectively, 
we have shown that 
$\Hom(\lambda_{\B}, \lambda_{\B})$
is an abelian algebra with dimension $k_1$ and it follows that 
$\lambda_{\B}$ decomposes into a direct sum of 
$k_1$ irreducible pieces, denoted by $\sigma
_i, 0\leq i\leq k_1-1.$\par
From $[\gamma \alpha_{\lambda_{\B}}] =
[\gamma \lambda\gamma]$ we have:
\[
\langle
\gamma \alpha_{\lambda_{\B}}, [g^i 
\lambda g^{-i}]\rangle = k_1, \quad 0\leq i\leq n_1-1.
\]
Note that
by (4) of Lemma \ref{key} we have 
\[ 
\langle
\gamma \alpha_{\lambda_{\B}}, g^i 
\lambda g^{-i}\rangle = \langle 
\alpha_{\lambda_{\B}}, g^i 
\lambda g^{-i}\rangle
.\] It follows that
$\alpha_{\lambda_{\B}}\succ k_1\bigoplus_{0\leq i\leq 
n_1-1}
[g^i  \lambda g^{-i}]$. On the other hand
\[
d({\alpha_{\lambda_{\B}}})= d({ 
\lambda_{\B}})= n d({\lambda}).
\] It follows that
$[\alpha_{\lambda_{\B}}]
= \bigoplus_{0\leq i\leq n_1-1}
k_1[g^i  \lambda g^{-i}]
.$
So we must have 
$\alpha_{\sigma_j} \succ g^i  \lambda g^{-i}$ for
some $i$ where $0\leq j\leq k_1-1.$  
By (3) of Lemma \ref{key}, $ [g\alpha_{\sigma_j}]=
[\alpha_{\sigma_j}g]$, so we must have 
$\alpha_{\sigma_j} \succ \bigoplus_{0\leq i\leq n_1-1}[g^i  
\lambda g^{-i}]$. In particular $d({\sigma_j})\geq 
n_1 d({\lambda})$. Since
$\sum_{0\leq j\leq k_1-1} d({\sigma_j})= k_1n_1 
d({\lambda})$,
it follows that 
\begin{equation}\label{al}
[\alpha_{\sigma_j}] = \bigoplus_{0\leq i\leq n_1-1}[g^i  
\lambda g^{-i}]
\end{equation}
and
$
\langle \alpha_{\sigma_0 }, \alpha_{\sigma_j} \rangle = n_1
.$
By (2) of Lemma \ref{key} we have
$\langle \alpha_{\sigma_0 }, \alpha_{\sigma_j} \rangle= \langle 
\sigma_0 \bar \sigma_j,
[1]+[\sigma]+...+[\sigma^{n-1}]\rangle =n_1$, it follows that the set 
$\{ \sigma_0, \sigma_1,..., \sigma_{k_1-1} \}$ is the same as 
$\{ \sigma_0, \sigma \sigma_0, ..., \sigma^{k_1-1} \sigma_0 \}$.
We will use $(\lambda; \sigma^i), 0\leq i \leq k_1-1$ 
to denote
$\sigma^i \sigma_0$ in the following. 
It follows from formula (\ref{al}) and (1) of Lemma \ref{key} that
$d({(\lambda; \sigma^i)})= n_1d(\lambda)$.
\par
The last part follows from  formula (\ref{al}) for $\alpha_{\sigma_j 
}$.
\endproof  
The following simple lemma will be used in  \S \ref{8.4} and \ref{8.5}.
\begin{lemma}\label{grading}
Let $\mu$ be an irreducible DHR representation of $\B$. Let
$i$ be any integer. Then:\par 
(1) $G(\mu,\sigma^i):=\e (\mu,\sigma^i) \e (\sigma^i,\mu) \in {\mathbb C},
$ 
$G(\mu,\sigma)^i=G(\mu,\sigma^i) $. Moreover $G(\mu,\sigma)^n=1$;\par
(2) If $\mu_1\prec \mu_2\mu_3$ with
$\mu_1,\mu_2,\mu_3$ irreducible, then $G(\mu_1,\sigma^i)=G(\mu_2, \sigma^i)
G(\mu_3,\sigma^i)$;\par
(3)   $\mu$ is untwisted if and only if  $G(\mu,\sigma)=1;$ \par
(4) $G(\bar \mu, \sigma^i)=\bar G(\mu, \sigma^i).$
\end{lemma}
\proof
We have  $G(\mu,\sigma^i)\in \Hom(\sigma^i \mu,\sigma^i \mu) \simeq  
{\mathbb C}$ since
$\sigma^i \mu$ is irreducible, and also 
$G(\mu,\sigma)^i= \e (\mu,\sigma^i) \e (\sigma^i,\mu)= G(\mu,\sigma^i)$, so
$G(\mu,\sigma)^n = 1$ since $[\sigma^n]=[1]$. 
If $\mu_1\prec \mu_2\mu_3$ with
$\mu_1,\mu_2,\mu_3$ irreducible, then $G(\mu_1, \sigma^i)=G(\mu_2,\sigma^i)
G(\mu_3, \sigma^i)$
by the  Braiding-Fusion equations (cf. \cite{R2}). 
For the second part, by Prop. \ref{nong}
$\mu$ is untwisted if and only if  $\alpha_\mu$ is a DHR representation
of $\D$. By the remark before Prop. 3.2 this is true if and only if 
$G(\mu, \sigma)=1$. The third part follows from $(2)$ and $G(1,\sigma^i)=1$.

\endproof 
Denote by $W$ the vector space whose basis consists of irreducible
components of  all $\alpha_\mu$ where $\mu$ are irreducible DHR 
representations of 
$\B,$   and 
$W_0$ (resp. $W_t$) the subspaces whose bases consist of  irreducible
components of   $\alpha_\mu$ where $\mu$ are irreducible non-twisted 
(resp. twisted) 
DHR  representations of 
$\B$ (relative to $\D$). 
The elements in the basis of $W_t$ are also called {\it twisted 
solitonic
sectors}. We note that  $W_0$ has a natural ring structure where the product
is the composition of sectors.  
Applying Prop. \ref{non} and Th.3.6 of \cite{BEK} we have
\begin{equation}\label{twist}
\dim W_t = \sum_{\lambda} (k_1-1)
\end{equation} 
So each 
$\lambda$ with nontrivial stabilizer contributes to
the twisted solitonic sectors. 
\subsection{Non-twisted sectors for the cyclic permutations}
For the rest of \S8, we will consider the case of cyclic permutations, 
i.e., we assume that  $\D:= \A\otimes\A...\otimes\A$ ($n$-fold tensor product) and 
$\B:=(\A\otimes \A...\otimes\A)^{\mathbb Z_n}$ the fixed point subnet 
of $\D$ under the action of cyclic permutations. Since we assume that $\D$
is completely rational, this is equivalent to assuming that $\A$ is
completely rational.  
We will denote by 
$(\lambda_1,...,\lambda_n)=\l_1\otimes\cdots\otimes\l_n$ the irreducible 
product representation of $\D$ associated with the irreducible 
representations $\l_1,\dots,\l_n$ of $\A$ and 
by $(\lambda_1,...,\lambda_n)_{\B}$ 
its restriction to $\B$. Note that $\mathbb Z_n$ acts on product sectors  $(\lambda_1,...,\lambda_n)$ 
naturally by cyclic permutations and
$[g (\lambda_1,...,\lambda_n)g^{-1})] 
=[(\lambda_{g(1)},...,\lambda_{g(n)})]. $
Assume that the stabilizer of $[(\lambda_1,...,\lambda_n)]$ is 
generated by
$g^{n_1}$ with $n_1k_1=n$. Then by Prop. \ref{nong} we have:
\begin{corollary}\label{non}
$(\lambda_1,...,\lambda_n)_{\B}$ decomposes into $k_1$ different 
irreducible pieces denoted by 
$(\lambda_1,...,\lambda_n; \sigma^i), 0\leq i\leq {k_1-1}$.
Moreover $ [\alpha_{(\lambda_1,...,\lambda_n; \sigma^i)}]=
\bigoplus_{0\leq k\leq {n_1-1}}[g^k(\lambda_1,...,\lambda_n)g^{-k} ]$, and
if  $[(\lambda_1,...,\lambda_n; \sigma^i)]= [(\mu_1,...,\mu_n; \sigma^j)]$
then  there exists an integer $l$ such that $\mu_{k}= 
\lambda_{g^l(k)}, 
1\leq k\leq n$.
\end{corollary}

\subsection{Topological twisted sectors for cyclic permutations}
\label{8.3}
Let us first determine the relevant ring structures of the
topological twisted sectors from Sect. \ref{tts}. 
Choose $\zeta$ to be the right boundary
point of the fixed interval $J_0$ in the anti-clockwise direction on the
circle. We can assume that $J_0$ is the interval $I$ as chosen in
Prop. \ref{non-v}.  Since $\pi_{f}^{(\z)}$ is a soliton, 
by the usual DHR argument \cite{DHR}, we can choose a soliton which is
unitarily equivalent to  $\pi_{f}^{(\z)}$ and restricts to an endomorphism
of $\D(J_0)$ (also cf. the paragraph before Prop. 4.5). We will denote this
endomorphism of $\D(J_0)$ by $\pi$. We note that $\tau_f$ is a DHR 
representation of $\B$ and we will denote by $\tau$ a DHR 
representation of $\B$ localized on the fixed interval $J_0$ which
is unitarily equivalent to $\tau_f$  and
the corresponding endomorphism of  $\B(J_0)$ obtained by
restriction to $\B(J_0)$. 
(Notations  differ here from the 
previously used ones: $\pi$ and $\t$ are sectors of factors). 

Note that by \cite{LR} we have
$[\tau]= [\gamma \pi\res \B(J_0)]$ as sectors of $\B(J_0)$. 
By $(d)$ of Th. \ref{indexformula} we have
$\langle\tau, \tau\rangle =n$, in fact 
$[\tau]= [\tau^{(0)}]+...+ [\tau^{(n-1)}]$.
So
$\langle  \gamma \pi\res \B(J_0),   \gamma \pi\res \B(J_0)\rangle=n$. 
As in the beginning of the proof of Prop. \ref{non}, we have
\begin{equation}\label{bd}
\langle  \gamma \pi\res \B(J_0),   \gamma \pi\res \B(J_0)\rangle
= \langle \pi,   \gamma \pi\gamma\rangle=n.
\end{equation}
By definition (cf.  \S \ref{3.1}) 
$[\gamma\alpha_{\tau}]=[\tau\gamma]=[\gamma\pi\gamma]$. We get
$\langle\gamma\alpha_{\tau}, \pi\rangle = 
\langle\gamma\pi\gamma , \pi\rangle=n.
$  
Since 
$[\gamma\alpha_{\tau}]=[\alpha_{\tau}]+...+[g^{n-1}\alpha_{\tau}]
$ and 
$\langle g^i\alpha_{\tau},  g^j\alpha_{\tau}\rangle = \delta_{ij}
\langle \alpha_{\tau},  \alpha_{\tau}\rangle ,
\forall 0\leq i,j\leq n-1$ (cf. (4) of Lemma \ref{key}), 
it follows that there exists an integer
$0\leq i\leq n-1$ such that 
$
\langle g^i\alpha_{\tau}, \pi\rangle = n$.
On the other hand since $d({\alpha_{\tau}})= d({\tau})=nd({\pi})$,
we must have $  [g^i\alpha_{\tau}]= n[\pi]$. Since 
$[\alpha_{\tau}]= 
[\alpha_{\tau^{(0)}}]+...+[\alpha_{\tau^{(n-1)}}]
$,
and $g^{-i}\pi$ is irreducible, 
we conclude that, for any $0\leq j\leq n-1$, we have 
$[\alpha_{\tau^{(j)}}] = [g^{-i} \pi]$.
 
Since $\alpha_{\tau^{(j)}}, \pi$ are solitons localized on
$J_0$ (cf. Prop. 3.1), using the next lemma we conclude that
\begin{equation}\label{at}
[\alpha_{\tau^{(j)}}] = [ \pi], \ 0\leq j\leq n-1.
\end{equation}
\begin{lemma}
Let $\pi_1, \pi_2$ be two solitons of $\D_0$ (The restriction of $\D$ to 
$S^1 \setminus \{ \z \},$ cf. \S3.0.1) localized on $J_0$. 
If 
$[\pi_1]= [g^{-i} \pi_2]$  as sectors of  $\D(J_0)$ for some integer $i$, then
$g^{-i}$ as a group element is the identity and 
 $[\pi_1]= [\pi_2]$  as sectors of $\D(J_0)$.
\end{lemma}
\begin{proof}
It is enough to prove that $g^{-i}$ as a group element is the identity.
Let $J_1 \subset J_0, J_1\neq J_0$ be an interval with $\z$ as a boundary
point. Let $J_2:= J_0\cap J_1'.$ Assume that $v$ is a unitary
in $\D(J_0)$ such that $\pi_1= \Ad v\cdot (g^{-i} \pi_2) $ on 
$\D(J_0)$. \par
Consider $\pi_1, \pi_2$ on $\D(J_1').$ Since $ \pi_1, \pi_2$ 
are solitons, and $\D(J_1')$ is a type III factor, we can find unitaries
$v_1,v_2$ such that on  $\D(J_1')$ we have
$\pi_1= \Ad{v_1},  \pi_2= \Ad{v_2}.$ Since
$  \pi_1, \pi_2$ are localized on $J_0$, it follows that
$v_1\in \D(J_0), v_2\in \D(J_0)$. \par
So on $ \D(J_2)$ we have 
$\Ad {v_1} = \Ad v\cdot \Ad{g^{-i}(v_2)} \cdot g^{-i}.$
Define $w:= g^{-i}(v_2)^* v^* v_1.$ Note that $w\in  \D(J_0)$, and
$ w x w^* = g^{-i}(x), \forall x\in \D(J_2).$ It follows that 
$w\in \D(J_0) \cap \B(J_2)'$. By (2) of Lemma 3.6 in \cite{Xu4} the
pair $\B\subset \D$ is strongly additive (cf. Definition 3.2 of 
\cite{Xu4} ) since we assume that $\D$
is strongly additive, and so  $\D(J_0') \vee \B(J_2)= \D(J_1')$ which
implies by Haag duality $\D(J_0) \cap \B(J_2)'= \D(J_1).$ Therefore
$w\in \D(J_1), g^{-i}(x)=x, \forall x\in \D(J_2),$ and so $g^{-i}$ as a
group element is the identity since one checks easily that the action of
the cyclic group on $\D (J_2)$ is faithful. 
\end{proof}
Note  that by Corollary \ref{canform} we have as 
sectors of $\D(J_0)$:
\begin{equation}\label{dispi}
[\bar\pi\pi]= \bigoplus_{\lambda_1,...,\lambda_n}\langle
\lambda_1\cdot\cdot\cdot\lambda_n, 1\rangle [(\lambda_1,...,\lambda_n)].
\end{equation}
where $\bar\pi$ is the conjugate sector of $\pi$.
From  $  [\alpha_{\tau^{(j)}}]= [\pi]$ we have
\begin{equation}\label{disa}
[ \alpha_{\bar\tau^{(j)}} 
\alpha_{\tau^{(j)}}
]= \bigoplus_{\lambda_1,...,\lambda_n}\langle
\lambda_1\cdot\cdot\cdot\lambda_n, 1\rangle [(\lambda_1,...,\lambda_n)].
\end{equation}
where we have also used $  [\alpha_{\bar\tau^{(j)}}] =[\bar
\alpha_{\tau^{(j)}}]$ (cf. (1) of Lemma \ref{key}).\par
Recall that the spins of $\tau^{(j)}$ are given in $(e)$ of 
Th. \ref{indexformula}. By (\ref{at}) 
$[\alpha_{\tau^{(j)}}]= [\alpha_{\tau^{(0)}}],$ by (2) of
Lemma \ref{key} we have 
$$
\sum_{0\leq l\leq n-1}\langle \tau^{(j)}, \sigma^l \tau^{(0)}\rangle
=1.
$$
Since both $\tau^{(j)}$ and $ \sigma^l \tau^{(0)}$ are irreducible,
we must have that
$
[\tau^{(j)}]=[ \sigma^{k(j)} \tau^{(0)}]$ 
where $k(\cdot)$ is a map from ${\mathbb Z}_n$ to itself.  $k(\cdot)$ is 
also one to one (hence onto) since if 
$[\sigma^{l_1}\tau^{(0)}]= [\tau^{(0)}]$ for some $0<l_1\leq n-1$, 
then by (2) of Lemma \ref{key} again
$
\langle \alpha_{\tau^{(0)}}, \alpha_{\tau^{(0)}}\rangle
= \sum_{ 0\leq l\leq n-1} \langle\tau^{(0)},  \sigma^l \tau^{(0)}
\rangle \geq 2
$
contradicting the fact that $\alpha_{\tau^{(0)}}$ is irreducible. 
We claim that in fact $k(j)= j k(1),  0\leq j\leq n-1.$ 
This follows essentially by the
grading Lemma \ref{grading}: by definition of the monodromy, 
$G(\sigma^{jk(1)} \tau^{(0)}, \sigma^{k(1)})=
G(  \tau^{(0)}, \sigma^{k(1)})$ 
because all $\sigma^j$'s have integer spins and are automorphisms.
From the monodromy equation (cf. \cite {R2}) we have 
\begin{equation}\label{k(1)}
G(\sigma^{jk(1)}   \tau^{(0)}, \sigma^{k(1)})= e^
{2\pi i({\rm spin}(\sigma^{(j+1)k(1) } \tau^{(0)}) - {\rm spin}
(\sigma^{jk(1)} 
\tau^{(0)}))}\ ,
\end{equation}
hence, modulo integers,
${\rm spin}(\sigma^{(j+1) k(1)} \tau^{(0)}) - {\rm spin}(\sigma^{jk(1)}  
\tau^{(0)})$ is
a constant independent of $0\leq j\leq n-1$. Since
$[\sigma^{k(1)}  \tau^{(0)}] =[\tau^{(1)}]
$,  ${\rm spin}(\sigma^{(j+1) k(1)} \tau^{(0)}) - {\rm spin}(\sigma^{jk(1)}  
\tau^{(0)})$ is equal to $\frac{1}{n}$ modulo integers. It follows that
${\rm spin}(\sigma^{j k(1)} \tau^{(0)})$ is  equal to the spin of 
$ \tau^{(j)}$ modulo integers. We conclude that
\[
[\sigma^{j k(1)} \tau^{(0)}] =   [\tau^{(j)}]\: \text{and}\:
jk(1)= k(j), \ 0\leq j\leq n-1\ .
\]
Since $k(\cdot)$ is one to one, the
greatest non-negative common divisor of $k(1)$ and $n$ must be 1.

In the following we define 
\[
G(\mu):= G(\mu, \sigma^{k(1)})
\] 
and will refer to $G(\mu)$ as {\it the grading } of $\mu$. 
Note that by definition $G(\tau^{(0)})= e^{\frac{2\pi i}{n}}.$ 
\par

Let $\l$ be a covariant representation of $\A$ and 
$\tau_{\lambda}= \pi_{\lambda}\res \B$ (cf. Prop. \ref{non-v}) the  
DHR representation of $\B$ obtained by restriction of 
$\pi_{\lambda}$. 
As in the beginning of this section, we denote
by $\pi_{\lambda}$ the endomorphism of $\D(J_0)$ 
obtained from the restriction
to  $\D(J_0)$ of a soliton unitarily equivalent to $ \pi_{\lambda}^{(\z)}$.
Note that an analogue of
$(d)$ of Th. \ref{indexformula} holds and $\tau_{\lambda}$  is a 
direct sum of $n$ DHR representations 
$\tau_{\lambda}^{(0)}, ...,\tau_{\lambda}^{(n-1)}$.

Note that $[\pi_{\lambda}]= [\pi\cdot(\lambda,1,...,1)]$ by 
Prop. \ref{non-v}, and it follows that 
$[\gamma\alpha_{\tau_{\lambda}}]=[ \gamma \pi 
(\lambda,1,...,1)\gamma]$.

By  the same argument as in the case when $\lambda=1$ above we 
have $[\alpha_{\tau_{\lambda}}]= n [g^{k} \pi (\lambda,1,...,1)]
=n [g^{k} \alpha_{\tau^{(0)}}\  (\lambda,1,...,1)]$
for some $0\leq k\leq n-1$, and by (3), (4) of Lemma \ref{key} again we 
have $k=0$ and 
\begin{equation}\label{alambda}
[\alpha_{\tau_{\lambda}^{(j)}}]= [\alpha_{\tau^{(0)}}\  
(\lambda,1,...,1)] =[(\lambda,1,...,1) \ \alpha_{\tau^{(0)}}]
=[\pi_{\lambda}],\ 0\leq j\leq n-1.
\end{equation}
From these equations we can prove the following:
\begin{theorem}\label{main}
$(1)$ $ [\pi_{\lambda}]= [\pi_{\mu}]$ as sectors of 
$\D(J_0)$ iff $\lambda\simeq \mu$ as DHR representations of $\A$;\par
$(2)$ $ [\tau_{\lambda}^{(l)}]=[\tau_{\mu}^{(j)}]$  iff 
$\lambda \simeq\mu$  as DHR representations of $\A$ and
$l=j$.
\end{theorem}
\proof
 $(1)$:
Since 
\[ 
[\pi_{\lambda}]=[\pi\cdot (\lambda,1,...,1)],\  
[\pi_{\mu}]=[\pi\cdot (\mu,1,...,1)]
\]
we have  $ [\pi_{\lambda}]= [\pi_{\mu}]$ iff 
$ [\pi (\lambda,1,...,1)]= [\pi (\mu,1,...,1)]$. It follows by 
Frobenius
duality and equation (\ref{dispi})
\[
\langle \bar\pi \pi, (\bar\lambda\mu,1,...,1) \rangle =1
= \langle 1, \bar\lambda\mu\rangle
.\]  
It follows that $[\lambda]=[\mu]$ as sectors of $\A(J_0)$. Since
$\A$ is strongly additive, it follows that   $\lambda\simeq \mu$ 
as DHR representations of $\A$. \par
 $(2)$:
It is sufficient to show that if  $ 
[\tau_{\lambda}^{(l)}]=[\tau_{\mu}^{(j)}] $ then  
$\lambda\simeq \mu$ as DHR representations of $\A$. 
Assume that  
$ [\tau_{\lambda}^{(l)}]=[\tau_{\mu}^{(j)}]. $
By  equation (\ref{alambda}) we have
\[
[\alpha_{\tau_{\lambda}^{(l)}}]= [\pi_{\lambda}],
\quad [\alpha_{\tau_{\mu}^{(j)}}]=  [\pi_{\mu}]
\]
and the proof follows from point $(1)$.
\endproof
We note that Th. \ref{main} is similar to the main theorems 
(Th. 3.9 and Th. 4.4)  of \cite{BDM} 
if one identifies
$\pi_{\lambda}$ with the twisted module in the sense of \cite{BDM}. 
Th. \ref{main}  supplies a class of twisted representation of the cyclic 
orbifold. 
In the next few sections we will show that these representations and 
their
variations give all the twisted representations in the case $n=2,3,4$.

\subsection{Case $n=2$}
\label{8.4}
When $n=2$, by (\ref{twist}) $\dim W_t$ is the same as the cardinality of 
the set $\{ \lambda  \}$. By (2) of Th. \ref{main} and (\ref{alambda})
$W_t$ has a basis $\{\alpha_{\tau_{\lambda}^{(0)}} \}$. 
If $\sigma_1$ is an irreducible twisted representation
of $\B$, it follows that 
$\alpha_{\sigma_1}=\bigoplus_\lambda C_\lambda 
\alpha_{\t_{\lambda}^{(0)}}$, where
$C_\lambda$ are positive integers. By eq. (\ref{k(1)}) and
 (2) of Lemma \ref{key} it follows that 
$\sigma_1$ must be some $\tau_{\lambda}^{(i)}$. One can 
also 
prove this by computing index of all known  DHR representations of 
$\B$  and check that 
they add up  to $\mu_{\B} = 4 \mu_{\D}$.  
Hence we have proved the
following:
\begin{proposition}\label{n=2}
When $n=2$ all the irreducible twisted representations of the fixed 
point 
net $\B$ are $\{\tau_{\lambda}^{(i)}\}$.
\end{proposition}
When $n=2$ we can determine completely the fusion rules of 
$\alpha_{\tau_{\lambda}^{(0)}}$ as follows: 
\begin{proposition}\label{fusion}
(1) $[\bar \alpha_{\tau_{\lambda}^{(0)}}]  =[ 
\alpha_{\tau_{{\bar\lambda}}^{(0)}}] ; $\par
(2) $[(\mu_1,\mu_2) \alpha_{\tau_{\lambda}^{(0)}}] = \bigoplus_\delta
\langle \mu_1\mu_2\lambda,
\delta\rangle  [\alpha_{\tau_{\delta^{(0)}}}]; $\par
(3)  $[\alpha_{\tau_{\lambda}^{(0)}}  \alpha_{\tau_{\mu}^{(0)}}]
=
\bigoplus_{\lambda_1,\lambda_2} \langle\lambda\mu\lambda_1, \lambda_2\rangle
[(\lambda_2,\bar\lambda_1)].$
\end{proposition} 
\proof
$(1)$: Note that $ [\alpha_{\tau_{\lambda}^{(0)}}  ]=[(\lambda,1)
\alpha_{\tau^{(0)}} 
]$, so it is sufficient to
show that $[\bar \alpha_{\tau^{(0)}}]  = 
[\alpha_{\tau^{(0)}}]     $. Here we give two different proofs.
Since $W_t$ is spanned by $\{ \alpha_{\tau_{{\bar\lambda}^{(0)}}}  
\}$, 
we must have that
$[\bar \alpha_{\tau^{(0)}}]= [  \alpha_{\tau_{\mu}^{(0)}} ]=
[(\mu,1)  \alpha_{\tau^{(0)}} ] $ for some $\mu$ (cf. 
(\ref{alambda})). 
From this
we have $d(\mu)=1$. So $\mu\lambda$ is irreducible for any $\lambda$. 
From 
\[
[\bar\alpha_{\tau_{\lambda}^{(0)}}] =[(\mu\bar\lambda,1) 
\alpha_{\tau^{(0)}} ],
[\bar\alpha_{\tau_{\lambda}^{(0)}}]= [\alpha_{\overline
{\tau_{\lambda}^{(0)}}}]
\]
we have 
\[
\langle
[\overline
{\tau^{(0)}_{\lambda}}], [\tau^{(0)}_{{\mu\bar\lambda}}]+
 [\tau^{(1)}_{{\mu\bar\lambda}}] \rangle =1
\]
and therefore
$\overline
{\tau^{(0)}_{\lambda}}$ is either $\tau^{(0)}_{{\mu\bar\lambda}}$ 
or
 $\tau^{(1)}_{{\mu\bar\lambda}}$. In any case the univalence
($=:\text{exp}(2\pi i\cdot\text{spin})$)  
$\omega_{{\tau_{\lambda}^{(0)}}}$ 
(cf. \cite{GL1}) of ${\tau_{\lambda}^{(0)}}$ must be the same as 
that of
$\tau_{\mu\bar\lambda}^{(0)}$ or
$\tau_{\mu\bar\lambda}^{(1)}$. 
Note that by (\ref{spin}) we have
\[
{\omega_{\tau_{\mu\bar\lambda}^{(0)}}}^2=
{\omega_{\tau_{\mu\bar\lambda}^{(1)}}}^2= \omega_{\mu\bar\lambda}
e^{\frac{2\pi ic}{16}}, 
{\omega_{\tau_{\lambda}^{(0)}}}^2= \omega_{\lambda}
e^{\frac{2\pi ic}{16}},  
\]
and therefore
$\omega_\lambda= \omega_{\mu\bar\lambda}, \forall \lambda$. It 
follows that
$\mu$ is degenerate (cf. \cite{R2}) and therefore $\mu$ is the 
vacuum representation since $\A$ is modular (cf. \cite{KLM}). 
This completes the first proof of 
 $[\bar \alpha_{\tau^{(0)}}]  = 
[\alpha_{\tau^{(0)}}]     .$
\par
For the second proof of  
 $[\bar \alpha_{\tau^{(0)}}]  = 
[\alpha_{\tau^{(0)}}] $, 
note that by (\ref{alambda}) 
$[\pi]= [\alpha_{\tau^{(0)}}]$.
By the remark after Prop. \ref{4.11}, we have 
$[\pi]= [\bar\pi]$.
So we have $[\alpha_{\tau^{(0)}}]= [\bar \alpha_{\tau^{(0)}}]$. 
\par
$(2)$: By (\ref{alambda}) 
we have
$[\alpha_{\tau_{\lambda}^{(0)}}]= [ (\lambda,1)\alpha_{\tau^{(0)}} 
].$ So 
\[ 
[(\mu_1,\mu_2) \alpha_{\tau_{\lambda}^{(0)}}] = \bigoplus_{\delta_1}
\langle \mu_1\lambda,\delta_1\rangle  [(\delta_1,\mu_2)
\alpha_{\tau^{(0)}} 
].
\]
Note that $ [(\delta_1,\mu_2)
\alpha_{\tau^{(0)}} 
]= [(\delta_1,1) (1,\mu_2)\alpha_{\tau^{(0)}} 
].$ We claim that $[(1,\mu_2)\alpha_{\tau^{(0)}} 
]= [(\mu_2,1)\alpha_{\tau^{(0)}} 
].$ In fact by (\ref{disa}) and Frobenius duality we have:
\begin{multline}
\langle (1,\mu_2)\alpha_{\tau^{(0)}}, (\mu_2,1)
\alpha_{\tau^{(0)}}\rangle
= \langle (\bar\mu_2,\mu_2), \bar\alpha_{\tau^{(0)}}\alpha_{\tau^{(0)}}\rangle =1 \\
\langle (1,\mu_2)\alpha_{\tau^{(0)}}, 
(1,\mu_2)\alpha_{\tau^{(0)}}\rangle
= \langle (1,\bar\mu_2)(1,\mu_2), \bar\alpha_{\tau^{(0)}}\alpha_{\tau^{(0)}}\rangle =1 \\
\langle (\mu_2,1)\alpha_{\tau^{(0)}}, (\mu_2,1)\alpha_{\tau^{(0)}}\rangle
= \langle (\bar\mu_2,1)(\mu_2,1), 
\bar\alpha_{\tau^{(0)}}\alpha_{\tau^{(0)}}\rangle =1
\end{multline}
It follows that  $[(1,\mu_2)\alpha_{\tau^{(0)}} 
]= [(\mu_2,1)\alpha_{\tau^{(0)}} 
].$ Hence 
$[(\delta_1,1) (1,\mu_2)\alpha_{\tau^{(0)}} 
]= [(\delta_1,1) (\mu_2,1)\alpha_{\tau^{(0)}} 
]= \bigoplus_{\delta} \langle\delta_1\mu_2, 
\delta\rangle[\alpha_{\tau_{\delta}^{(0)}} 
].$ So we have 
\[
[(\mu_1,\mu_2) \alpha_{\tau_{\lambda}^{(0)}}] = \bigoplus_{\delta_1,\delta}
\langle \mu_1\lambda,\delta_1\rangle  
\langle \delta_1\mu_2,\delta\rangle 
[\alpha_{\tau_{\delta}^{(0)}}
]= \bigoplus_{\delta}
\langle \mu_1\mu_2\lambda,\delta\rangle 
[\alpha_{\tau_{\delta}^{(0)}} 
].
\] \par
$(3)$: We have
\[
[\alpha_{\tau_{\lambda}^{(0)}}  \alpha_{\tau_{\mu}^{(0)}}]
=[(\lambda,1)(\mu,1) 
\alpha_{\tau^{(0)}}  \alpha_{\tau^{(0)}}]
=\bigoplus_{\lambda_1} [(\lambda,1)(\mu,1) (\lambda_1,\bar\lambda_1)]
=\bigoplus_{\lambda_1,\lambda_2} \langle\lambda\mu\lambda_1, \lambda_2\rangle
[(\lambda_2,\bar\lambda_1)]
\]
where we have used  (\ref{alambda}) in the first equality, 
the first part of the proposition and (\ref{disa}) in the second 
equality. 
\endproof
\par
Before concluding this subsection, we note that $\pi_{\mu}$ can be 
defined also for a  reducible sector $\mu$ of $\A$ and we clearly have
\[
\pi_{\mu} = \bigoplus_{\delta}\langle \mu,\delta\rangle\pi_{\delta}\ ,
\]
where $\delta$ runs on the irreducible sectors of $\A$.

Hence Proposition \ref{fusion} can be equivalently formulated, with 
the notations in Sect. \ref{tts}, as follow:
\begin{itemize}
\item[$(1)$] $\bar\pi_{\l}\simeq \pi_{\bar\l}$,
\item[$(2)$] $(\mu_1\otimes\mu_2)\cdot\pi_\l
\simeq \pi_{\mu_1 \mu_2 \l}$,
\item[$(3)$] $\pi_{\l}\pi_{\mu}\simeq
\bigoplus_{\delta}\l\delta\otimes\mu\bar\delta$\ ,
\end{itemize}
where $\l,\mu,\mu_1 ,\mu_2$ and $\delta$ are irreducible.

(1) is proved in Prop. \ref{pi}, (2) follows from the equality 
$(\mu\otimes\iota)\cdot\pi_\l = (\iota\otimes\mu)\cdot\pi_\l=\pi_{\mu\l}$
and (3) follows by Cor. \ref{canform}. Note that the composition of 
two twisted solitons is a DHR sector.
\subsection{Case $n=3$} 
\label{8.5}
By (\ref{twist}) when $n=3$ $\dim W_t$ is twice the cardinality of  the 
set $\{ \lambda \}$. 
We claim that in this case (unlike the case $n=2$)
$[\alpha_{\tau_{\lambda}^{(0)}}] \neq [\overline{
\alpha_{\tau_{\mu}^{(0)}}}]$. 
If not, by Frobenius duality, (\ref{alambda}) and (2) of Lemma \ref{key} 
we have
\[
\langle (\bar\lambda \bar\mu,1,1), \alpha_{\tau^{(0)}}^2\rangle
= \langle \gamma(\bar\lambda \bar\mu,1,1)\res \B(J_0), 
{\tau^{(0)}}^2\rangle
=1.
\]
It follows that 
${\tau^{(0)}}^2$ contains some untwisted DHR representation of 
$\B$.   
Note that $G(\tau^{(0)})^2= e^{\frac{2\pi i}{2}}=-1$, so by 
Lemma  
\ref{grading}  we have arrived at a contradiction. 
Hence by counting we conclude that
$W_t$ is spanned by $ \{ \alpha_{\tau_{\lambda}^{(i)}}, 
\alpha_{\bar\tau_{\lambda}^{(i)}} \},$ and by the same argument 
as in the proof 
of Prop. \ref{n=2} we have:
\begin{proposition}
All the irreducible twisted representations of $\B$ in the case $n=3$ 
are
$\tau_{\lambda}^{(i)}$ and $\bar \tau_{\lambda}^{(i)}, 0\leq 
i\leq 2.$
\end{proposition}
\subsection{Case $n=4$}
By (\ref{twist}) in this case $\dim W_t= |\{ 
(\lambda_1,\lambda_2,\lambda_1,\lambda_2), 
\lambda_1\neq \lambda_2, \}| + 3  |\{ 
(\lambda,\lambda,\lambda,\lambda) \}
|$. One question is how to construct additional sectors 
corresponding to
$(\lambda_1,\lambda_2,\lambda_1,\lambda_2)$. We notice that there is
an intermediate fixed point net $\C$ between $\B$ and $\D$ such that
$\C$ is the fixed point subnet of $\B$ under the action of $g^2$. In 
fact
$\C$ is  fixed point subnet of 
$\D=(\A\otimes \A)\otimes (\A\otimes \A)$ under the natural cyclic
${\mathbb Z_2}$ action. So we can apply the results of 
 \S \ref{8.4} to the pair $\C\subset \D$. 
Now the representations of $\A\otimes \A$ are
labeled by $(\lambda_1,\lambda_2)$, and so we label the solitons for the
pair $\C\subset \D$ by $\pi_{(\lambda_1,\lambda_2)}$ and its 
restriction
to $\C$ (a DHR representation of $\C$) 
by $\tau_{(\lambda_1,\lambda_2)}$. Recall from \S \ref{8.3}
that
 $\tau_{(\lambda_1,\lambda_2)}$ is a direct sum of two 
irreducible DHR
representations denoted by  $\tau^{(0)}_{(\lambda_1,\lambda_2)}$ 
and 
 $\tau^{(1)}_{(\lambda_1,\lambda_2)}$. We will denote by 
 $\tau^{(i)}_{(\lambda_1,\lambda_2), \B}$ the DHR representations 
of $\B$ obtained by restricting   
$\tau^{(i)}_{(\lambda_1,\lambda_2)}$
to $\B$, $i=0,1$. Note that $\C$ is invariant under the automorphism 
induced
by cyclic permutation $g$ and the $\B$ is the fixed point subnet 
under this
action. Applying Prop. \ref{nong} to $\B\subset \C$ we have
\begin{equation}\label{e1}
[\alpha^{\B\uparrow \C}_{\tau^{(0)}_{(\lambda_1,\lambda_2),\B}}] 
= [\tau^{(0)}_{(\lambda_1,\lambda_2)}]+
[ g \tau^{(0)}_{(\lambda_1,\lambda_2)} g^{-1}]
\end{equation}
where $\B\uparrow C$ indicates the induction from $\B$ to $\C$ (note 
that an horizontal arrow has been used in \cite{Xu4}). 
By Lemma 3.3 of \cite{Xu3} we have
\begin{equation}\label{e2}
[\alpha^{\B\uparrow \D}_{\tau^{(0)}_{(\lambda_1,\lambda_2),\B}}] 
= [\alpha^{\C\uparrow 
\D}_{\tau^{(0)}_{(\lambda_1,\lambda_2)}}]+
[ \alpha^{\C\uparrow \D}_{g \tau^{(0)}_{(\lambda_1,\lambda_2)} 
g^{-1}}]
\end{equation}
By $(3)$ of Lemma \ref{key}  as sectors  
$\alpha^{\B\uparrow \D}_{\tau^{(0)}_{(\lambda_1,\lambda_2), 
\B}}$ 
commutes with $g$ since $[g]$ is a 
subsector of the canonical endomorphism $\gamma$ from $\D$ to $\B$.  
So we must have
\begin{equation}\label{e3}
[g\alpha^{\C\uparrow 
\D}_{\tau^{(0)}_{(\lambda_1,\lambda_2)}}g^{-1}] 
= [\alpha^{\C\uparrow 
\D}_{\tau^{(0)}_{(\lambda_1,\lambda_2)}}] \; 
{\mathrm {or}}\;
[ \alpha^{\C\uparrow \D}_{g \tau^{(0)}_{(\lambda_1,\lambda_2)} 
g^{-1}}]
\end{equation} 
As in the proof of (\ref{bd})  and using (\ref{alambda}) we have
\begin{equation}\label{e4}
\langle   \tau_{{(\lambda_1,\lambda_2)},\B}    ,  
\tau_{{(\lambda_1,\lambda_2)}, \B}\rangle
=\langle\alpha^{\C\uparrow 
\D}_{\tau^{(0)}_{(\lambda_1,\lambda_2)}},  
\gamma 
\alpha^{\C\uparrow \D}_{\tau^{(0)}_{(\lambda_1,\lambda_2)}}
\gamma \rangle
\end{equation}
By using (\ref{e3}), (\ref{e4})  we conclude that 
$\tau_{{(\lambda_1,\lambda_2), \B}}$ 
is a direct sum of four distinct irreducible
pieces iff 
\[
[g\alpha^{\C\uparrow 
\D}_{\tau^{(0)}_{(\lambda_1,\lambda_2)}}g^{-1}] 
= [\alpha^{\C\uparrow \D}_{\tau^{(0)}_{(\lambda_1,\lambda_2)}}]
,\] and
a direct sum of two distinct irreducible pieces iff 
\[
[g\alpha^{\C\uparrow 
\D}_{\tau^{(0)}_{(\lambda_1,\lambda_2)}}g^{-1}] 
= 
[ \alpha^{\C\uparrow \D}_{g 
\tau^{(0)}_{(\lambda_1,\lambda_2)}g^{-1}}
]\neq  
[\alpha^{\C\uparrow \D}_{\tau^{(0)}_{(\lambda_1,\lambda_2)}}]
.\]
On the other hand, applying Prop. \ref{nong} to the pair
$\B\subset \C$, we know that $\tau_{(\lambda_1,\lambda_2), \B}$
is a direct sum of four irreducible pieces iff 
$[g \tau^{(i)}_{(\lambda_1,\lambda_2),\B} g^{-1}] =
 [\tau^{(i)}_{(\lambda_1,\lambda_2),\B}], i=0,1,$
and a direct sum  of two distinct irreducible pieces iff 
$[g \tau^{(i)}_{(\lambda_1,\lambda_2),\B} g^{-1}] \neq
 [\tau^{(i)}_{(\lambda_1,\lambda_2),\B}], i=0,1,$ and
$[g \tau^{(0)}_{(\lambda_1,\lambda_2),\B} g^{-1}] \neq
 [\tau^{(1)}_{(\lambda_1,\lambda_2),\B}]
$. So we have that 
$
[g\alpha^{\C\uparrow 
\D}_{\tau^{(0)}_{(\lambda_1,\lambda_2)}}g^{-1}] 
= [\alpha^{\C\uparrow \D}_{\tau^{(0)}_{(\lambda_1,\lambda_2)}}]
$ iff 
$[g \tau^{(i)}_{(\lambda_1,\lambda_2),\B} g^{-1}] =[
 \tau^{(i)}_{(\lambda_1,\lambda_2),\B}], i=0,1$, and
$[ \alpha^{\C\uparrow \D}_{g 
\tau^{(0)}_{(\lambda_1,\lambda_2)}g^{-1}} 
] = [g\alpha^{\C\uparrow 
\D}_{\tau^{(0)}_{(\lambda_1,\lambda_2)}}g^{-1}]. $
In particular 
\begin{multline}\label{e5}
[ \alpha^{\C\uparrow \D}_{g \tau^{(0)}_{(\lambda_1,\lambda_2)}g^{-1}}]
=[g\alpha^{\C\uparrow\D}_{\tau^{(0)}_{(\lambda_1,\lambda_2)}g^{-1}}]
= [g\alpha^{\C\uparrow \D}_{\tau^{(0)}_{(1,1)}}g^{-1} g 
(\lambda_1,
\lambda_2,1,1) g^{-1}] \\
= [\alpha^{\C\uparrow \D}_{\tau^{(0)}_{(1,1)} g^{-1}}
(\lambda_2, 1,1,\lambda_1)]
\end{multline}
Note that $[\alpha^{\C\uparrow \D}_{g 
\tau^{(0)}_{(1,1)}g^{-1}}]= [g\alpha^{\C\uparrow 
\D}_{\tau^{(0)}_{(1,1)}}g^{-1}]$, and so 
$g\tau^{(0)}_{(1,1)} g^{-1}$ is a twisted DHR representation  of $\C$ 
(relevant to $\D$). Applying Prop. \ref{n=2} to the pair
$\C\subset\D$ we have
\begin{equation}\label{e6}
[\alpha^{\C\uparrow \D}_{\tau^{(0)}_{(1,1)}}g^{-1}
]=
[\alpha^{\C\uparrow \D}_{\tau^{(0)}_{(1,1)}}
(\sigma_1,\sigma_2,1,1)] 
\end{equation}
for some $(\sigma_1,\sigma_2)$. By (\ref{e5}) we have
$[ \alpha^{\C\uparrow \D}_{g \tau^{(0)}_{(\lambda_1,\lambda_2)} 
g^{-1}}]
= [ \alpha^{\C\uparrow \D}_{\tau^{(0)}_{(1,1)}} 
(\sigma_1\lambda_2, \sigma_2,1,
\lambda_1)]$, and by $(2)$ of Prop. \ref{fusion} we have
$ [ \alpha^{\C\uparrow \D}_{\tau^{(0)}_{(1,1)}}
(\sigma_1\lambda_2, \sigma_2,1,
\lambda_1)]
=  [ \alpha^{\C\uparrow 
\D}_{\tau^{(0)}_{(\sigma_1\lambda_2 ,\sigma_2\lambda_1 )}} 
]$. Hence
\begin{equation}\label{e7}
[ \alpha^{\C\uparrow \D}_{g \tau^{(0)}_{(\lambda_1,\lambda_2)} 
g^{-1}}]
=  [ \alpha^{\C\uparrow 
\D}_{\tau^{(0)}_{(\sigma_1\lambda_2 ,\sigma_2\lambda_1 )}} 
]
\end{equation}
By $(2)$ of Lemma \ref{key} 
we have that
$g \tau^{(0)}_{(\lambda_1,\lambda_2)} 
g^{-1} \simeq \tau^{(i)}_{(\sigma_1\lambda_2 ,\sigma_2\lambda_1 )},
$ where $i=0$ or $i=1$, as DHR representations of $\C$. Notice that
$\omega_{g \tau^{(0)}_{(\lambda_1,\lambda_2)} g^{-1}}=
\omega_{\tau^{(0)}_{(\lambda_1,\lambda_2)}}$ which can be checked
directly from the definition of univalence (cf. \cite{GL2}). 
Alternatively one 
can prove this as follows. First if $g \tau^{(0)}_{(\lambda_1,\lambda_2)} 
g^{-1} \simeq \tau^{(0)}_{(\lambda_1,\lambda_2)}$ then we have 
nothing to prove. If 
 $[g \tau^{(0)}_{(\lambda_1,\lambda_2)} 
g^{-1}] \neq [\tau^{(0)}_{(\lambda_1,\lambda_2)}]$, applying 
Prop. \ref{non} to the pair $\B\subset \C$ we know that 
$g \tau^{(0)}_{(\lambda_1,\lambda_2)} 
g^{-1}$ and $\tau^{(0)}_{(\lambda_1,\lambda_2)}$ restricts to the
same DHR representation of $\B$, and so they 
must have the same univalence by
Lemma 6.1 of \cite{BE2}.\par
So we have  
\begin{equation}\label{e8}
\omega_{g \tau^{(0)}_{(\lambda_1,\lambda_2)} g^{-1}}=
\omega_{\tau^{(0)}_{(\lambda_1,\lambda_2)}}=
\omega_{\tau^{(i)}_{(\sigma_1\lambda_1,\sigma_2\lambda_2)}}  
\end{equation}
where $i=0$ or $1$.
As in the first proof of 
(1) of Prop. \ref{fusion}, from (\ref{e8})
we have
$\omega_{(\sigma_1\lambda_2, \sigma_2\lambda_1)}=\omega_{(\lambda_1, 
\lambda_2)}, \forall (\lambda_1,\lambda_2).$
It follows that $(\sigma_1,\sigma_2) $ is degenerate. 
Therefore
$ (\sigma_1,\sigma_2)=(1,1),$ and 
\[
[ \alpha^{\C\uparrow \D}_{\tau^{(0)}_{(1,1)}} ]= 
[ g\alpha^{\C\uparrow \D}_{\tau^{(0)}_{(1,1)}} g^{-1}]
.\]
By (\ref{e7}) we have
\[
[ \alpha^{\C\uparrow \D}_{\tau^{(0)}_{(\lambda_2,\lambda_1)}} ]= 
[ \alpha^{\C\uparrow \D}_{g\tau^{(0)}_{(\lambda_1,\lambda_2)} g^{-1}}]
,\]
and by (\ref{alambda}) we have
\begin{equation}\label{e9}
[ \alpha^{\C\uparrow \D}_{\tau^{(i)}_{(\lambda_2,\lambda_1)}} ]= 
[ \alpha^{\C\uparrow \D}_{g\tau^{(i)}_{(\lambda_1,\lambda_2)} 
g^{-1}}], \;  i=0,1.
\end{equation}
If $\lambda_1=\lambda_2$, by (\ref{e9}) and the remark before (\ref{e5})
we must have $g\tau^{(i)}_{(\lambda_1,\lambda_2)} g^{-1}\simeq 
\tau^{(i)}_{(\lambda_1,\lambda_2)}, i=0,1$. 
Apply Prop. \ref{nong} to
the pair $\B\subset\C$
we know that $\tau^{(i)}_{(\lambda_1,\lambda_1),\B}$ is a direct sum of 
two 
distinct irreducible pieces denoted by 
$\tau^{(i,j)}_{(\lambda_1,\lambda_1),\B}, i,j=0,1.
$
\par
If 
$\lambda_1\neq \lambda_2$, then from  (\ref{e9}) and (2) of Lemma \ref{key}
we have that
$g\tau^{(i)}_{(\lambda_1,\lambda_2)} g^{-1}\simeq 
\tau^{(j)}_{(\lambda_2,\lambda_1)}
$ where $0\leq j\leq 1$, and
$[g\tau^{(i)}_{(\lambda_1,\lambda_2)} g^{-1}]\neq 
[\tau^{(i)}_{(\lambda_1,\lambda_2)}]$.
We may choose our labeling so that 
$g\tau^{(i)}_{(\lambda_1,\lambda_2)} g^{-1}\simeq 
\tau^{(i)}_{(\lambda_2,\lambda_1)}
.$
Apply Prop. \ref{nong} to
the pair $\B\subset\C$
we know that $\tau^{(i)}_{(\lambda_1,\lambda_1),\B}$ is an irreducible
DHR representation of $\B$, and
$\tau^{(i)}_{(\lambda_1,\lambda_2),\B}$ are isomorphic to
$\tau^{(i)}_{(\lambda_2,\lambda_1),\B}$ as DHR representations of 
$\B, i=0,1.$
By definitions 
we have $G(\tau^{(i)}_{(\lambda_1,\lambda_2),\B})^2= 
G(\tau^{(i,j)}_{(\lambda_1,\lambda_1),\B})^2=1,$ since 
$\tau^{(i)}_{(\lambda_1,\lambda_2),\B}$ and $
\tau^{(i,j)}_{(\lambda_1,\lambda_1),\B}$
are  non-twisted representations of $\B$ relevant to $\C$, and
$\B$ is the fixed point subnet of $\C$ under the $\mathbb Z_2$ action.
So
these representations  are  different from 
$\tau_{\lambda}^{(i)}$ whose grading is $e^{\frac{2\pi i }{4}}$
or $\bar \tau_{\lambda}^{(i)}$ whose grading is $e^{\frac{6\pi i 
}{4}}$. \par
We note that by applying Prop. \ref{nong} to $\B\subset \C$ we have 
\[
d({\tau^{(i,j)}_{(\lambda_1,\lambda_1),\B}})= 
d({\tau^{(i)}_{(\lambda_1,\lambda_1)}}),\
d({\tau^{(i)}_{(\lambda_1,\lambda_2),\B}})= 
2 d({\tau^{(i)}_{(\lambda_1,\lambda_2)}}), \quad \lambda_1\neq \lambda_2.
\]
Applying  (\ref{alambda}) and (c) of Th. \ref{indexformula} to 
$\C\subset \D$ we have
\[
d^2({\tau^{(i)}_{(\lambda_1,\lambda_1)}})= 4 d^2({(\lambda_1,\lambda_1)})
\mu_{\A}^2, \
d^2({\tau^{(i)}_{(\lambda_1,\lambda_2)}})= 4 d^2({(\lambda_1,\lambda_1)})
\mu_{\A}^2. 
\]
Hence we know the indices of these known twisted representations 
$\tau^{(i,j)}_{(\lambda_1,\lambda_1),\B}, \tau^{(i)}_{(\lambda_1,\lambda_2),\B}$
of $\B$ 
(relevant to $\D$). By Prop. \ref{non} we also know
the indices of non-twisted  representations of $\B$ 
relevant to $\D$. 
One can check easily that the sum of these indices add up to $\mu_{\B}
=16 \mu_{\D}= 16  \mu_{\A}^4$. 
By \cite{KLM} we have therefore identified all the
irreducible DHR representations of $\B$. In particular we have proved
the following:
\begin{proposition}
All the irreducible twisted DHR 
representations of $\B$ (relevant to $\D$) are 
$\tau_{\lambda}^{(i)},$ 
$\bar \tau_{\lambda}^{(i)}, 0\leq i\leq 3, $ 
$\tau^{(i,j)}_{(\lambda,\lambda),\B}, 0\leq i,j\leq 1,$
and
$\tau^{(i)}_{(\lambda_1,\lambda_2), \B}, \lambda_1\neq \lambda_2, 
0\leq i\leq 1,$ where as DHR representations 
$\tau^{(i)}_{(\lambda_1,\lambda_2), \B}$ are isomorphic to
$\tau^{(i)}_{(\lambda_2,\lambda_1),\B}.$
\end{proposition}
We note that our construction of 
$\tau^{(i)}_{(\lambda_1,\lambda_2), \B}$ 
and  
$\tau^{(i,j)}_{(\lambda,\lambda), \B}$ 
can be generalized to non-prime $n$ case. 
\subsection{Comments on the case of a general $n$}
To motivate our discussion let us first consider the case when $\A$ 
is holomorphic,
i.e. when $\mu_\A=1$. In this case $\D$ is also holomorphic, and $\D$ 
has
only one irreducible representation (the vacuum) labeled by 
$(1,...,1)$.
In this case $\dim W_t = n-1$. Note that $\alpha_{\tau^{(0)}}\in 
W_t$ is a periodic automorphism, and we let 
$k\geq 1$ be the least integer such that  
$[\alpha_{\tau^{(0)}}^k]= [1].$
By Lemma \ref{grading} we must have $n|k$. On the other hand we must 
have
$k\leq n$ since $\dim W=n$. So we conclude that $k=n,$
and $W$ is spanned by 
$\{ \alpha_{\tau^{(0)}}^{i}, 0 \leq i \leq n-1 \}$ and all the 
irreducible representations of $\B$ are given by $\sigma^j  
{\tau^{(0)}}^{i},
0\leq i, j\leq n-1.$ So in the holomorphic case all twisted 
representations
of $\B$ are generated by $ \tau^{(0)}$ and $\sigma$ via fusion. This
example shows that it is an interesting question to determine the 
nature of ``composed'' sectors $\alpha_{\tau^{(0)}}^k$ 
($k\in {\mathbb N}$) in the general
case as we have done for the case $n=2$ in \S \ref{8.3}.\par
For general completely rational $\A$, we note that the grading
$G(\tau^{(0)})= 
e^{\frac{2\pi i}{n}}$ by the remark after the definition of grading 
in \S8.3. 
Now if $\sigma_1$ is an irreducible twisted DHR representation of $\B$,
by Lemma \ref{grading} the grading $G(\sigma_1)$ is a complex number 
such that
$G(\sigma_1)^n=1$. Assume that 
$G(\sigma_1)=  e^{\frac{-2\pi ki}{n}}, 1\leq k\leq n-1.$ Let $\sigma_2$ 
be any irreducible DHR representation of  $\B$ such that 
${\tau^{(0)}}^k\succ
\sigma_2$. By $(1)$ of Lemma \ref{grading} $G(\sigma_2)=  e^{\frac{2\pi 
ki}{n}}$
and if $\mu\prec \sigma_1 \sigma_2$ is an irreducible DHR 
representation  of $\B$, 
then $G(\mu)=1.$  It follows from Lemma \ref{grading} that 
$\mu$ is non-twisted whose 
nature
is determined in Cor. \ref{non}. By using Frobenius duality, we conclude 
that
$\bar\sigma_1\prec \bar\mu{\tau^{(0)}}^k$. 
This observation shows once again the
importance of  $\tau^{(0)}$ and suggests that it is an interesting
question to determine the nature of 
${\tau^{(0)}}^k$ ($k\in {\mathbb N}$) in general
case. This question is related to the question in the previous 
paragraph
by Lemma \ref{key}.
\section{Generalizations and the case of two-dimensional nets}
Results and proofs in this paper remain valid with weaker assumptions. 
We replace axiom $\bf D$ by the following ones:
\begin{itemize}
\item{\it Reeh-Schlieder property}: $\Om$ is cyclic for $\A(I)$,
$I\in\I$.
\item{\it Modular PCT}: The modular conjugation of $(\A(S^+),\Omega)$ 
corresponds to the reflection $z\mapsto\bar z$ of $S^1$. (By 
M\"{o}bius covariance the modular conjugations associated with all 
intervals have then a geometric meaning.)
\item{\it Factoriality}: $\A(I)$ is a factor for all $I\in\I$.
\item{\it Equivalence between local and global intertwiners}: If 
$\mu,\nu$ are finite-index endomorphisms localized in the interval 
$I$, then $\text{Hom}(\mu,\nu)=\text{Hom}(\mu_I,\nu_I)$ as in 
\cite{GL2}. 
\end{itemize}
If $\C$ is a local conformal net on the two-dimensional Minkowski 
spacetime $\mathbb R^2$ (see \cite{KL2}), 
let $\A$ be the restriction of $\C$ to the 
time-zero axis: $\A(I)\equiv\C(\O)$ where $\O$ is the double cone 
with basis $I$. Then $\A$ satisfies all the above properties hence 
our results do apply. In particular we then have:
\begin{theorem}
If $\C$ is a local conformal net on the two-dimensional Minkowski 
spacetime. The following are equivalent:
\begin{itemize}
\item[$(i)$] $\A$ is not completely rational;
\item [$(ii)$] $\sum_i d(\r_i)=\infty$ (sum over all irreducible 
sectors);
\item [$(iii)$] $(\A\otimes\A)^{\rm flip}$ has an irreducible sector 
with infinite dimension.
\end{itemize}
\end{theorem}
The rest of our results have analogous extensions.

\appendix
\section{$\Mob^{(n)}$ covariance in the strongly additive case}
\label{covar}
Let $E$ be a symmetric $n$-interval of $S^1$, namely 
$E\equiv\sqrt[n]{I}$ for some $I\in\I$. With $I_0 , I_1, \cdots 
I_{n-1}$ 
the $n$ connected component of $E$, by the split property we have 
a natural isomorphism
\[
\chi_E: \A(I_0)\otimes\A(I_1)\otimes\cdots\otimes\A(I_{n-1})
\to\A(I_0)\vee\A(I_1)\vee\cdots\vee\A(I_{n-1}) = \A(E) \ .
\]
A state of the form
\[
\f\equiv (\f_0\otimes\f_1\otimes\cdots\otimes\f_{n-1})\cdot\chi^{-1}_E
\]
on $\A(E)$, where $\f_k$ is a normal faithful state on $\A(I_k)$ and 
$\f_k=\f_0\cdot \Ad(U(R(2k\pi/n))$, is called  a \emph{rotation 
invariant product state}. 

We state here a formula for the modular group of 
$\A(E)$, that extends to the general case the formula by
Schroer and Wiesbrock \cite{SW} in the example of the 
$U(1)$-current algebra, see  \cite{L7}.
 
\begin{proposition}
\label{modgroup}
There is a rotation 
invariant product state $\f$ on $\A(E)$ such that the corresponding 
modular group $\s^{\f}$ 
of $\A(E)$ is given by
\[
\s_t^{\f}=\Ad U^{(n)}(\Lambda_I(-2\pi t))\res\A(E)
\]
where $\Lambda_I$ is the the lift to \Mob$^{(n)}$ of one parameter 
subgroup of 
{\rm \Mob} of 
generalized dilation associated with $I$ (see \cite{GL2}) and 
$U^{(n)}=U\cdot M^{(n)}$ is the unitary representation of 
\Mob$^{(n)}$.
\end{proposition}
\begin{corollary}
Let $\A$ be a strongly additive
local conformal net on $S^1$ with the split 
property. Then every representation of $\A$ with finite index
is {\rm \Mob}$^{(n)}$-covariant with positive energy, for all $n\in\mathbb N$.
\end{corollary}
\proof
As $\A$ is strongly additive, every finite index sector is 
$\Mob$-covariant with positive energy by \cite{GL1}.

Fix $n$ and let $\Ad U^{(n)}(g)$ be the action of $\Mob^{(n)}$ 
on $\A$ given in Sect. \ref{first}. 
Let $\r$ be a finite-index localized endomorphism. We may
assume $\r$ to be localized in an interval which is a connected component 
of a symmetric $n$-interval $E=\sqrt[n]{I}$.

With $\{\Lambda_I^{(n)}(t)\}\subset\Mob^{(n)}$ the 
one-parameter dilation subgroup, 
denote by $\a_t\equiv \Ad U^{(n)}\cdot\Lambda_I^{(n)}(-2\pi t)$ the 
corresponding 
rescaled action on $\A$. 

We have to show that $\r_t\equiv\a_t\cdot\r\cdot\a^{-1}_t$ is 
equivalent to $\r$
for every $t\in\mathbb R$, namely that there is a unitary $z_t\in\A(E)$ 
such that
\begin{equation}\label{cov}
\r =\Ad z_t\cdot\a_t\cdot\r\cdot\a^{-1}_t\ ;
\end{equation}
having the covariance with respect to 
$\Lambda_I^{(n)}$, by changing the interval $I$ we then get the
covariance with respect to $\Mob^{(n)}$.

By the Prop. \ref{modgroup} $\a$ restricts to the modular 
automorphism group of $\A(E)$ with respect to $\f$. With $\Phi_{\r}$ 
the left inverse of $\r\restriction \A(E)$, by \cite{L3} the
Connes \cite{Con} cocycle $z_t=(D\phi\cdot\Phi_\r :\f)_t\in\A(E)$ 
satisfies 
\[
\r(x) =\Ad z_t\cdot\a_t\cdot\r\cdot\a^{-1}_t (x), \quad x\in\A(E),
\]
hence we obtain eq. (\ref{cov}) by strong additivity.
\endproof
\section{Frobenius reciprocity for global intertwiners}
In this section we show that Th. 3.21 of \cite{BE} (also cf. (4)
of Lemma \ref{key}) holds for global intertwiners when $\N$ is a 
conformal subnet of conformal net $\M$ with finite index. Note that
we do not assume strong additivity conditions for the net $\N$ as in 
\cite{BE},  but 
we consider global intertwiners.  We will use the
notations in \S 3 of \cite{BE} and refer the reader to \cite{BE} for
unexplained notations. Fix an interval $I_0.$  Let $\M$ be a 
conformal net on a Hilbert space $\H$ and 
$\lambda_1, \lambda_2$ be two DHR representations of 
$\M$ localized on $I_0$. Define 
\[
\Hom (\lambda_1, \lambda_2):= 
\{ x\in B(\H)| x \lambda_{1,J}(m) =\lambda_{2,J}(m) x, \forall m\in 
\M(J), \forall J \}. 
\]
$\Hom (\lambda_1, \lambda_2)$ will be called the space of 
{\it global intertwiners} from $\lambda_1$ to $\lambda_2.$ Its 
dimension
will be denoted by $\langle \lambda_1,  \lambda_2\rangle.$
The elements of $\Hom (\lambda_{1,I_0}, \lambda_{2,I_0}):= \{ x\in \M(I_0)|
 x \lambda_{1,I_0}(m) =\lambda_{2,I_0}(m) x, \forall m\in \M(I_0) \}$ 
are referred to as {\it local intertwiners}  from $\lambda_1$ to 
$\lambda_2$
(localized on $I_0$). Note that by Haag duality one obviously has
$\Hom (\lambda_1, \lambda_2)\subset \Hom (\lambda_{1,I_0}, 
\lambda_{2,I_0}).$
The following simple lemma tells us when a local intertwiner is global.
\begin{lemma} 
Let $I$ be an open interval which contains the closure of $I_0$. If 
$
x\in \M (I_0) \cap
\Hom (\lambda_{1,I}, \lambda_{2,I}),
$
then $x\in \Hom (\lambda_1, \lambda_2)$.
\end{lemma}
\proof
By definition we can cover any interval $J$ by $I$ and a finite 
number of
intervals $I_k$ such that $I_k\in I_0'.$ By the additivity of 
$\M$ we have $\M(J)\subset \M(I) \vee(\vee_k \M(I_k))$ and the lemma follows 
from the definitions.
\endproof
Now let $\lambda, \beta$ be DHR representations of $\N$ and $\M$ 
respectively localized in $I_0,$ and $\sigma_\beta$ be the DHR 
representation
of $\N$ localized on $I_0$ obtained from restriction of $\beta$ to 
$\N$. Assume that $\alpha_\lambda$ is a DHR representation of $\M$. 
We
have the following theorem:
\begin{theorem}
$
\langle \alpha_\lambda, \beta\rangle_\M = 
\langle \lambda, \sigma_\beta\rangle_\N. 
$
\end{theorem}
\proof We will adapt the proof of Th.  3.21 of \cite{BE}.  Choose an 
interval $I$ as in the Lemma B.1.  We can choose a $Q$-system 
$(\g_I,v,w)$ for the inclusion $\N(I)\subset\M(I)$ so that $\g_{\tilde 
I}$ extends to a canonical endomorphism of $\M(\tilde I)$ into $
\N(\tilde I)$ for all intervals $\tilde I\supset I$ so that 
$(\g_{\tilde I},v,w)$ $Q$-system for $\N(\tilde I)\subset\M(\tilde I)$.

First we show the inequality ``$\leq$''.  Let $t\in \Hom ( 
\alpha_\lambda, \beta)$.  By Haag duality we have $t\in \M(I_0)$ and 
$r=\gamma (t) w \in \N(I_0)$.  The argument on Page 25 of \cite{BE} 
shows that $r\in \Hom ( \lambda_I, (\sigma_\beta)_I)$.  By Lemma B.1 
we have $r\in \Hom ( \lambda, \sigma_\beta)$.  By Lemma 3.4 of 
\cite{BE} the map $t\rightarrow r$ is injective, thus ``$\leq$'' is 
proved.  \par 
We now turn to prove  ``$\geq$''.  Suppose that $r\in \Hom ( 
\lambda, \sigma_\beta)$ is given.  By Haag duality $r\in \N(I_0)$, and 
so $t=v^*r\in \M(I_0), s=\gamma(t)\in \N(I_0).  $ Clearly $s\in \Hom 
((\theta\lambda)_I, (\sigma_\beta)_I)$ since $r$ is a global 
intertwiner.  It follows by Lemma B.1 that $s$ is also a global 
intertwiner, and so Lemma 3.20 of \cite{BE} applies.  The rest of the 
proof is exactly the same as the proof on Page 26 of \cite{BE}.  
\endproof
    
\noindent{\bf Acknowledgments.}
The first named author would like to thank S. Carpi, F. Fidaleo, Y. Kawahigashi and
L. Zsido for comments. He also thanks Sorin Popa for the invitation 
and warm hospitality at UCLA in May 2003 while this work was in progress.    

\end{document}